\documentclass[11pt]{amsart}

\usepackage[centertags]{amsmath}
\usepackage{pstricks}
\usepackage{pst-node}
\usepackage{hyperref}

\numberwithin{equation}{section}

\renewcommand{\MR}[1]{\href{http://www.ams.org/mathscinet-getitem?mr=#1}{MR #1}}
\newcommand{\arXiv}[1]{\href{http://front.math.ucdavis.edu/#1}{arXiv:#1}}

\newcommand{\bi}{\bigskip}
\newcommand{\sm}{\smallskip}

\newcommand{\ee}{\end{equation}}
\newcommand{\eea}{\end{eqnarray}}
\newcommand{\bean}{\begin{eqnarray*}}
\newcommand{\eean}{\end{eqnarray*}}

\newif\ifpctex\newcommand{\noi}{\noindent}

  \newtheorem{theorem}{Theorem}
  \newtheorem{definition}{Definition}[section]

  \newtheorem{cond}[definition]{Condition}
  \newtheorem{proposition}[definition]{Proposition}
  \newtheorem{lemma}[definition]{ Lemma}
  \newtheorem{cor}[definition]{Corollary}
  \newtheorem{ass}[definition]{Assumption}
  
  \newcommand{\beCond}[2]{\Rand{\vspace{0,6cm}\tt #1}\begin{cond}[#2]\label{#1}}

\newcommand{\R}{\mathbb{R}}  

\newcommand{\N}{\mathbb{N}}

\newcommand{\be}[1]{\begin{equation}\label{#1}}
\newcommand{\bew}[1]{\begin{equation*}\label{#1}}
\newcommand{\bea}[1]{\begin{eqnarray}\label{#1}}

\newcommand{\beL}[2]{\begin{lemma}[#2]\label{#1}}
\newcommand{\beD}[2]{\begin{definition}[#2]\label{#1}}
\newcommand{\beT}[2]{\begin{theorem}[#2]\label{#1}}
\newcommand{\beP}[2]{\begin{proposition}[#2]\label{#1}}
\newcommand{\beC}[2]{\begin{cor}[#2]\label{#1}}

\newcommand{\length}{\mbox{\rm length}}
\newcommand{\lengths}{\mbox{\rm lengths}}
\newcommand{\shape}{\mbox{\rm shape}}
\newcommand{\re}[1]{(\ref{#1})}
  
\begin{document}

\title[real trees]{Rayleigh processes, real trees,
and root growth with re-grafting}

\author{Steven N.\ Evans}
\address{Department of Statistics \#3860 \\
  University of California at Berkeley \\
367 Evans Hall \\
Berkeley, CA 94720-3860 \\
U.S.A}
\email{evans@stat.Berkeley.EDU}
\thanks{SNE supported in part by NSF grants DMS-0071468 and DMS-0405778,
and a Miller Institute
for Basic Research in Science research professorship}

\author{Jim Pitman}
\address{Department of Statistics \#3860 \\
  University of California at Berkeley \\
367 Evans Hall \\
Berkeley, CA 94720-3860 \\
U.S.A}
\email{pitman@stat.Berkeley.EDU}
\thanks{JP supported in part by NSF grants  DMS-0071448 and DMS-0405779}

\author{Anita Winter}
\address{Mathematisches Institut \\
Universit\"at Erlangen--N\"urnberg \\
Bismarckstra\ss{}e $1\frac12$ \\
91054 Erlangen \\
GERMANY}
\email{winter@mi.uni-erlangen.de}
\thanks{AW supported by a DFG Forschungsstipendium}

\date{\today}



\begin{abstract}
The real trees form a class of metric spaces that extends
the class of trees with edge lengths 
by allowing behavior such as infinite
total edge length and vertices with infinite branching degree.
Aldous's Brownian continuum random tree, the random tree-like object
naturally associated with a standard Brownian excursion, may
be thought of as a random compact real tree.  The continuum random
tree is a scaling limit as $N \rightarrow \infty$
of both a critical Galton-Watson tree conditioned to have total
population size $N$ as well as a uniform random rooted combinatorial tree
with $N$ vertices.  The Aldous--Broder algorithm is a Markov chain
on the space of rooted combinatorial trees with $N$ vertices that
has the uniform tree as its stationary distribution.  We construct
and study a Markov process on the space of all rooted compact real trees
that has the continuum random tree as its stationary distribution
and arises as the scaling limit as $N \rightarrow \infty$
of the Aldous--Broder chain.  A key technical ingredient in this work
is the use of a pointed Gromov--Hausdorff distance to metrize
the space of rooted compact real trees.\\
\\
\\
Berkeley Statistics Technical Report No. 654 (February 2004),
revised October 2004. To appear in {\em Probability Theory and Related Fields}.
\end{abstract}
\maketitle

\newpage

\tableofcontents

\newpage
\section{Motivation and background}

It is shown in \cite{Ald91a, MR93f:60010, Ald93} (see also
\cite{MR2001g:60211, csp})
that a suitably rescaled family of Galton-Watson trees, conditioned to 
have total population size $n$, converges as $n\to\infty$ 
to the {\em Brownian continuum random tree (CRT)}, which can be
thought of as the tree inside a {\em standard Brownian excursion}
(more precisely, twice a standard Brownian excursion
if one follows Aldous's choice of re-scaling). 
Aldous describes a procedure for representing trees as
closed subsets of $\ell^1$, and
convergence is here in the sense of weak convergence of probability
measures on the this space of closed subsets
equipped with the Hausdorff distance.

The Brownian CRT can be obtained as an almost
sure limit by growing finite trees (that is,
trees with finitely many leaves and finite total branch
length) in continuous
time as follows (at all times $t \ge 0$ the procedure will
produce a rooted tree ${\mathcal R}_t$ with total edge length $t$):
\begin{itemize}
\item 
Write $\tau_1, \tau_2, \ldots$ for the successive arrival times
of an inhomogeneous Poisson process with arrival rate $t$ at time $t \ge 0$.
Call $\tau_n$ the {\em $n^{\mathrm{th}}$ cut time}.
\item
Start at time $0$
with the $1$-tree (that is a line segment with two ends), 
${\mathcal R}_0$, of length zero 
(${\mathcal R}_0$ is ``really'' the  {\em trivial tree}
that consists of one point only, but thinking this way helps
visualize the dynamics more clearly for this semi-formal
description). Identify one end of ${\mathcal R}_0$ as the root.
\item
Let this line segment grow at unit speed until the first cut time $\tau_1$.  
\item
At time $\tau_1$ pick a point uniformly on the segment
that has been grown so far.  Call this point the
{\em first cut point}.
\item
Between time $\tau_1$ and time $\tau_2$, 
evolve a tree with $3$ ends
by letting  a new branch growing away from the 
first cut point at unit speed.
\item
Proceed inductively:  Given the $n$-tree 
(that is, a tree with $n+1$ ends), 
${\mathcal R}_{\tau_n-}$, pick the {\em $n$-th cut point} uniformly
on ${\mathcal R}_{\tau_n-}$ to give an $n+1$-tree,
${\mathcal R}_{\tau_n}$, with one edge of length zero,
and for $t\in[\tau_n,\tau_{n+1}[$, let
${\mathcal R}_t$ be the tree obtained from ${\mathcal R}_{\tau_n}$ by
letting a branch grow away from the $n^{\mathrm{th}}$ cut
point with unit speed. 
\end{itemize}
 The tree
${\mathcal R}_{\tau_n-}$ has the same distribution as the
subtree of the CRT that arises from sampling twice the
standard Brownian excursion at $n$ i.i.d. uniform points
on the unit interval, and 
the Brownian CRT is the limit (with respect to a suitable notion of
convergence) of the increasing family of rooted finite trees
$({\mathcal R}_{t})_{t \ge 0}$.

Again using the cut times 
$\{\tau_1,\tau_2,\ldots\}$, a closely related 
way of growing rooted trees $({\mathcal T}_t)_{t \ge 0}$
can be described as follows (here again the tree at time $t$
will have total edge length $t$): 
\begin{itemize}
\item 
Start with the $1$-tree (with one end identified as the {\em root}
and the other as a {\em leaf}), ${\mathcal T}_0$, of length zero. 
\item
Let this segment grow at unit speed on the time interval $[0,\tau_1[$, 
and for $t\in[0,\tau_1[$ let ${\mathcal T}_t$
be the rooted $1$-tree that has its points labeled by the interval
$[0,t]$ in such a way that the root is $t$ and the leaf is $0$. 
\item
At time $\tau_1$ sample 
the {\em first cut point} uniformly along the tree ${\mathcal T}_{\tau_1-}$, 
prune off the piece of  ${\mathcal T}_{\tau_1-}$ that
is above the cut point (that is, prune off the interval of points
that are further away from the root $t$ than the first cut point).
\item 
Re-graft the pruned segment such that its cut end  and the root are
glued together. Just as we thought of ${\mathcal T}_0$ as a tree with two points,
(a leaf and a root) connected by an edge of length zero, we take 
${\mathcal T}_{\tau_1}$ to be the the rooted $2$-tree obtained by
``ramifying'' the root ${\mathcal T}_{\tau_1-}$  into two points (one of
which we keep as the root) that are joined by an edge of length zero.
\item 
Proceed inductively: Given the labeled and rooted $n$-tree, 
${\mathcal T}_{\tau_{n-1}}$, for $t\in[\tau_{n-1},\tau_n[$, let ${\mathcal T}_t$ be
obtained by letting the edge containing the root
grow at unit speed so that the points in ${\mathcal T}_t$ correspond
to the points in the interval $[0,t]$ with $t$ as the root. 
At time $\tau_{n}$, the $n^{\mathrm{th}}$ cut point is sampled randomly 
along the edges of the $n$-tree, ${\mathcal T}_{\tau_n-}$, 
and the subtree above the cut point (that is the subtree of points
further away from the root than the cut point)
is pruned off and re-grafted so that
its cut end and the root are glued together. The root is then ``ramified''
as above to give an edge of length zero leading from the root to the
rest of the tree.
\end{itemize}

\setlength{\unitlength}{0.7pt}
\begin{picture}(350,345)(-20,-235)
\thinlines
\put(0,30){\circle*{5}}
\put(100,30){\line(0,1){50}}
\put(15,25){\mbox{${\mathcal R}_0$}}
\put(115,25){\mbox{${\mathcal R}_{\tau_1-}$}}
\put(215,25){\mbox{${\mathcal R}_{\tau_1}$}}
\put(310,25){\mbox{${\mathcal R}_{\tau_2}-$}}

\put(410,25){\mbox{${\mathcal R}_{\tau_2}$}}
\put(200,30){\line(0,1){50}}
\put(97,30){\line(1,0){6}}
\put(197,30){\line(1,0){6}}
\put(197,80){\line(1,0){6}}
\put(292,30){\line(1,0){6}}
\put(200,60){\circle*{5}}
\put(295,30){\line(0,1){50}}
\put(295,60){\line(-1,0){25}}
\put(292,80){\line(1,0){6}}
%
\put(395,30){\line(0,1){50}}
\put(395,50){\circle*{5}}
\put(395,60){\line(-1,0){25}}
\put(392,80){\line(1,0){6}}
\put(370,57){\line(0,1){6}}
\put(392,30){\line(1,0){6}}
%
%
\put(15,-75){\mbox{${\mathcal T}_0$}}
\put(115,-75){\mbox{${\mathcal T}_{\tau_1}-$}}
\put(215,-75){\mbox{${\mathcal T}_{\tau_1}$}}
\put(320,-75){\mbox{${\mathcal T}_{\tau_2}-$}}
\put(420,-75){\mbox{${\mathcal T}_{\tau_2}$}}
\put(0,-70){\circle*{5}}
\put(100,-70){\line(0,1){50}}
\put(100,-40){\circle{3}}
\put(97,-20){\line(1,0){6}}
%
\put(200,-70){\line(0,1){30}}
\put(200,-70){\line(1,1){15}}
\put(212,-55){\line(1,0){6}}
\put(197,-40){\line(1,0){6}}
\put(197,-70){\line(1,0){6}}
%
\put(300,-45){\line(1,1){15}}
\put(297,-15){\line(1,0){6}}
\put(312,-29){\line(1,0){6}}
\put(300,-70){\line(0,1){55}}
\put(300,-35){\circle{3}}
%
\put(400,-45){\line(1,1){15}}
\put(397,-35){\line(1,0){6}}
\put(412,-29){\line(1,0){6}}
\put(400,-70){\line(0,1){35}}
\put(397,-70){\line(1,0){6}}
\put(400,-70){\line(-1,1){15}}
\put(382,-54){\line(1,0){6}}

\put(-20,-120){\makebox{
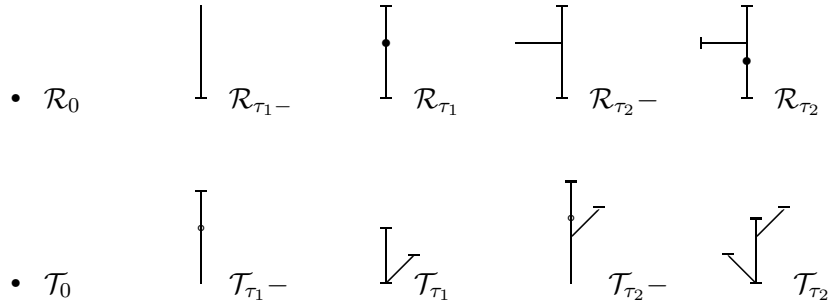
\begin{minipage}[t]{12cm}{\em {\sc Figure~1} illustrates how the 
 tree-valued processes $({\mathcal R}_t;\,t\ge 0)$ and
 $({\mathcal T}_t;\,t\ge
 0)$ evolve. (The bold dots re-present an edge of length zero, while the 
 small dots indicate the position of the cut point that is going to show 
 up at the next moment.)} 
\end{minipage}}}
\end{picture}

The link between these two dynamics for growing trees is
provided by Proposition \ref{P:Aldousconnect}, where we show that 
${\mathcal T}_{\tau_n-}$  has
the same law as ${\mathcal R}_{\tau_n-}$ for each $n$.

The process $({\mathcal T}_t)_{t \ge 0}$ is clearly a time-homogeneous
Markov process and we can run its dynamics (which we will
refer to as {\em root growth with re-grafting}) starting with
any finite tree.  The resulting process evolves
via alternating deterministic root growth
and random jumps due to re-grafting and is an
example of a {\em piecewise-deterministic Markov process}. 
A general framework for such processes was
introduced in \cite{MR87g:60062}  as an
abstraction of numerous examples in queueing and control theory,
and this line of research was extensively developed in the
subsequent monograph \cite{MR96b:90002}. 
A more general formulation in terms of
martingales and additive functionals can be found in \cite{JacSko96}.
Some other appearances of such processes are 
\cite{MR99j:60109,MR2002h:60179,MR2000g:60125,MR94e:60062,MR91d:60169}.
We note also that pruning and re-grafting operations such as the one
we consider here play an important role in algorithms that attempt
to reconstruct optimal phylogenetic trees from data by moving through
tree space as part of a hill-climbing or simulated annealing procedure
(see, for example, \cite{Fel03}).

The crucial feature of the root growth with re-grafting
dynamics is that they have a simple projective structure:
If one follows the evolution of the points in a rooted subtree of the
initial tree along with that of the points added at later times
due to root growth, then these points together form a rooted subtree at each
period in time and this subtree evolves autonomously according
to the root growth with re-grafting dynamics.

The presence of this projective structure 
suggests that one can make sense of the
notion of running the root growth with re-grafting
dynamics starting from an initial ``tree'' that has exotic behavior
such as infinitely many leaves, points with infinite branching,
and infinite total edge length -- provided that this ``tree'' can
be written as the increasing limit of a sequence of finite
trees in some appropriate sense.  Moreover, by the remarks
above about the relationship between the processes ${\mathcal R}$
and ${\mathcal T}$, this extended process should have a stationary
distribution that is related to the Brownian CRT, and the stationary
distribution should be the limiting distribution for the extended
process starting from any initial state. 

One of our main objectives is to give rigorous statements and proofs
of these and related facts.

Once the extended process has been constructed, 
we gain a new perspective on objects such as standard
Brownian excursion and the associated random triangulation of the circle 
(see \cite{MR95i:60007, MR95b:60011, Ald00}). 
For example, suppose we follow the height (that is, distance from the
root) of some point in the initial
tree.  It is clear that this height evolves autonomously as a one-dimensional
piecewise-deterministic Markov process that:
\begin{itemize}
\item 
increases linearly at unit speed (due to growth at the root), 
\item makes jumps at rate $x$ when it is in state $x$ 
(due to cut points falling
on the path that connects the root to the point we are following),
\item
jumps from state $x$ to a point
that is uniformly distributed on $[0,x]$ (due to re-grafting at the root).
\end{itemize}
We call such a process a {\em Rayleigh process} because,
as we will show in Section~\ref{S:Rayleigh}, this process converges to the 
{\em standard Rayleigh} stationary distribution $\mathbf{R}$ on $\R_+$
given by
\begin{equation*}
\mathbf{R}(]x,\infty[) = e^{-x^2/2}, \quad  x \ge 0,
\end{equation*}
(thus $\mathbf{R}$ is also the distribution of the Euclidean
length of a two-dimensional
standard Gaussian random vector or, up to a scaling constant,
the distribution of the distance to the closest point to the origin
in a standard planar Poisson process). Now, if 
$B^{\mathrm ex}:=\{B^{\mathrm ex}_u;\,u\in [0,1]\}$ is standard Brownian excursion
and $U$ is an independent uniform
random variable on $[0,1]$, then there is a
valid sense in which $2B^{\mathrm ex}_{U}$ has the law of the height 
of a randomly sampled
leaf of the Brownian CRT, and this accords with the well-known result
\be{height}
   {\bf P}\{2B^{\mathrm ex}_{U}\in dx\}=\mathbf{R}(dx).
\ee


In order to extend the root growth with re-grafting dynamics to
infinite trees we will need to fix on a suitable class
of infinite trees and a means of measuring distances
between them. 
Our path to extending the definition of a tree to accommodate
the ``exotic'' behaviors mentioned above will be the one
followed in the so-called
{\em T-theory} (see \cite{Dre84, MR97e:05069, Ter97}). T-theory 
takes finite trees to be just metric spaces
with certain characteristic properties and then
defines a more general class of tree-like metric spaces called
{\em real-trees} or {\em $\R$-trees}.  We note that one of the
primary impetuses for the development of T-theory was
to provide mathematical tools for concrete problems
in the reconstruction of phylogenies.  We also note that
$\R$-trees have been objects of intensive study in geometric
group theory (see, for example, the surveys
\cite{MR89d:57012, MR92e:20017, MR94e:57020, MR2003b:20040}
and the recent book \cite{MR2003e:20029}).
Diffusions on an $\R$-tree were investigated in \cite{MR2002e:60128}.
Some of the results on the space of $\R$-trees
obtained in this paper have already been found useful in the study
scaling limits of Galton--Watson branching processes
in \cite{DuqLeG04}.

Once we have an extended notion of trees
as just particular abstract metric spaces (or, more correctly, 
isometry classes of metric spaces), we need a means of assigning
a distance between two metric spaces, and this is provided
by the {\em Gromov-Hausdorff distance}.  This distance originated
in geometry as a means of making sense of intuitive notions
such as the convergence to Euclidean space 
of a re-scaled integer lattice as the grid size approaches zero 
or the convergence to Euclidean space of a sphere when viewed 
from a fixed point (for example, the North Pole)
as the radius approaches infinity.  
Our approach is thus rather different to
Aldous's in which trees are viewed as closed subsets
of $\ell^1$ via a particular choice of embedding and distances are
measured using the familiar Hausdorff metric.  Although
our use of the Gromov-Hausdorff distance to metrize the space
of $\R$-trees turns out to quite elegant and easy to work with,
it also appears to be rather novel, and thus much of our work in this
paper is directed towards establishing facts about the
structure of this space.  However, we think that the
resulting mathematics is interesting in its own right
and potentially useful in other investigations where trees
have hitherto been coded as other objects such as paths (for
example, \cite{MR2001g:60211, MR1954248}).
We remark in passing that the papers \cite{MR90k:57015, MR90d:57015}
are an application of the Gromov-Hausdorff distance to the study of
$\R$-trees that is quite different to ours.

We note that there is quite a large literature 
on other approaches to ``geometrizing'' and
``coordinatizing'' spaces of trees. 
The first construction of codes 
for labeled trees without edge-length 
goes back to 1918:  Pr\"ufer \cite{Pru18}
sets up a bijection between labeled trees of size $n$ 
and the points of
$\{1,2,\ldots,n\}^{n-2}$. 
Phylogenetic trees are identified with points 
in matching polytopes in \cite{DiaHol98}, and \cite{BilHolVog01} 
equips the space of finite phylogenetic trees with a fixed number of leaves
with a metric that makes it a cell-complex
with non-positive curvature. 

The plan of the rest of the paper is as follows. We collect some
results on the set of isometry classes of rooted compact $\R$-trees
and the properties of the Gromov-Hausdorff distance in Section~\ref{rgp}.
We construct the extended root growth with re-grafting process in
Section~\ref{S:construct} via a procedure that is roughly analogous
to building a discontinuous Markov process in Euclidean space as the solution
of a stochastic differential equation with respect to a sufficiently
rich Poisson noise.  This approach is particularly well-suited to
establishing the strong Markov property.  In Section~\ref{S:Aldousconnect}
we establish the fact claimed above that 
${\mathcal T}_{\tau_n-}$  has
the same law as ${\mathcal R}_{\tau_n-}$ for each $n$.
We prove in Section~\ref{S:recurstat}
that the extended root growth with re-grafting process is recurrent
and convergent to the continuum random tree stationary distribution.
We verify that the extended process has a Feller semigroup in
Section~\ref{S:Feller}, and show in Section~\ref{S:Aldous-Broder} 
that it is a re-scaling limit of the Markov chain appearing in the
Aldous--Broder algorithm for simulating a uniform rooted tree on
some finite number of vertices.  We devote Section~\ref{S:Rayleigh} to
a discussion of the Rayleigh process described above.

\section{$\R$-trees}
\label{rgp}

\subsection{Unrooted trees}
\label{unroot}
A complete metric space $(X,d)$ is said to be an {\em $\R$-tree} if it
satisfies the following axioms:\bi

\noi{\bf Axiom~1 (Unique geodesics) } For all 
$x,y\in X$ there exists a unique isometric embedding
  $\phi_{x,y}:[0,d(x,y)]\to X$ such that $\phi_{x,y}(0)=x$
and $\phi_{x,y}(d(x,y))=y$.\bi

\noi
{\bf Axiom~2 (Loop-free) } For every injective continuous map 
$\psi:[0,1]\to X$ one
has $\psi([0,1])=\phi_{\psi(0),\psi(1)}([0,d(\psi(0),\psi(1))])$.\bi

We refer the reader to 
(\cite{Dre84, DreTer96, MR97e:05069, Ter97}) for 
background on $\R$-trees.  A particularly useful fact is
that a metric space $(X,d)$ is an $\R$-tree if and only if
it is complete, path-connected, and satisfies the so-called 
{\em four point condition}, that is, 
\be{four}
\begin{aligned}
   d(x_1,&x_2)+d(x_3,x_4)
  \\
 &\le
   \max\{d(x_1,x_3)+d(x_2,x_4), \; d(x_1,x_4)+d(x_2,x_3)\}
\end{aligned}
\ee
for all $x_1,\ldots,x_4\in X$.  


Recall that the {\em Hausdorff distance}
between two subsets $A_1$, $A_2$ of a metric
space $(X,d)$ is defined as
\be{haus}
   d_{\mathrm H}(A_1,A_2)
 :=
   \inf\{\varepsilon>0;\,A_1\subseteq U_\varepsilon(A_2)\mbox{ and }
   A_2\subseteq U_\varepsilon(A_1)\},
\ee
where 
\be{epsnbhd}
U_{\epsilon}(A):=\{x\in X;\,d(x,A)\le\varepsilon\}. 
\ee
Based on this notion of distance between
closed sets, we define the {\em Gromov-Hausdorff distance},
$d_{{\mathrm{GH}}}(X_1,X_2)$, between
two metric spaces $(X_1,d_{X_1})$ and $(X_2,d_{X_2})$ as the 
infimum of $d_{\mathrm H}(X_1^\prime,X_2^\prime)$ over all metric spaces
$X_1^\prime$ and $X_2^\prime$ that are isomorphic to $X_1$ and $X_2$,
respectively,
and that are subspaces of some common metric space $Z$
(compare \cite{MR2000d:53065, MR2000k:53038, BurBurIva01}).
The   Gromov-Hausdorff distance defines a finite
metric on the space of all isometry classes of compact metric spaces
(see, for example, \ Theorem~7.3.30 in \cite{BurBurIva01}). \sm

Let $({\bf T},d_{{\mathrm{GH}}})$ be the metric space of isometry 
classes of compact 
$\R$-trees equipped with $d_{{\mathrm{GH}}}$. 
We will elaborate ${\bf T}$ slightly to incorporate the
notion of rooted trees, and this latter space
of rooted trees will be the state space of the 
Markov process having root growth with re-grafting
dynamics that we are going to construct. 
We will be a little loose and
sometimes refer to an $\R$-tree as an element of ${\bf T}$ rather
than as a class representative of an element.\bi

\noi
{\em Remark }\quad 
As we remarked in the Introduction, Aldous's approach to formalizing
the intuitive notion of an infinite tree and putting a
metric structure on the resulting class of objects is to work
with particular closed subsets of  
$\ell^1$ and to measure distances using the Hausdorff 
metric on closed sets.  Seen in the light of our approach,
Aldous's approach uses the distance between
two particular representative elements for 
the isometry classes of a pair of trees rather than the two that minimize 
the Hausdorff distance. In general this leads to greater
distances and hence a topology that is stronger than ours.\bi

The following results says that, at the very least, ${\bf T}$
equipped with the Gromov-Hausdorff distance is a ``reasonable''
space on which to do probability theory.

\begin{theorem}\label{T0} The metric space
$({\bf T},d_{{\mathrm{GH}}})$ is 
complete and separable.
\end{theorem}\sm

Before we prove Theorem~\ref{T0}, we point out that 
a direct application
of the stated definition of the Gromov-Hausdorff distance
requires an optimal embedding into a new metric space $Z$.   
While this definition is conceptually appealing and
builds on the more familiar Hausdorff distance between sets,
it turns out to often not be so useful for explicit
computations in concrete examples.
A re-formulation of the
Gromov-Hausdorff distance is suggested by the following
observation.  Suppose that two spaces $(X_1,d_{X_1})$ and $(X_2,d_{X_2})$
are close in the Gromov-Hausdorff distance as witnessed by isometric
embeddings $f_1$ and $f_2$ into some common space $Z$.  The map that
associates each point in $x_1 \in X_1$ to a point in $x_2 \in X_2$ such
that $d_Z(f_1(x_1),f_2(x_2))$ is minimal should then be close
to an isometry onto its image, and a similar remark holds
with the roles of $X_1$ and $X_2$ reversed.

In order to quantify the observation of the previous
paragraph, we require some more notation. 
A subset ${\Re}\subseteq X_1\times X_2$
is said to be a {\em correspondence} between 
sets $X_1$ and $X_2$ if for each $x_1\in X_1$ there exists at least one 
$x_2\in X_2$ such that $(x_1,x_2)\in{\Re}$, and for each 
$y_2\in X_2$ there exists at least one $y_1\in X_1$ such that
$(y_1,y_2)\in{\Re}$. Given metrics $d_{X_1}$ and $d_{X_2}$ on
$X_1$ and $X_2$, respectively,  the {\em distortion} of ${\Re}$ is
defined by 
\be{dis}
   {\mathrm{dis}}({\Re})
 :=
   \sup\{|d_{X_1}(x_1,y_1)-d_{X_2}(x_2,y_2)|;\,(x_1,x_2),
   (y_1,y_2)\in{\Re}\}.
\ee
Then
\be{GH}
   d_{{\mathrm{GH}}}((X_1,d_{X_1}),(X_2,d_{X_2}))
 =
   \frac{1}{2}\inf_{{\Re}}{\mathrm{dis}}({\Re}),
\ee
where the infimum is taken over all correspondences ${\Re}$
between $X_1$ and $X_2$ (see, for example,
 Theorem~7.3.25 in \cite{BurBurIva01}).\bi

The following result is also useful in the proof of Theorem~\ref{T0}.

\begin{lemma}
\label{L:Tclosed}
The set ${\bf T}$ of compact $\R$-trees 
is a closed subset of the space of compact
metric spaces equipped with the Gromov-Hausdorff distance.
\end{lemma}

\begin{proof}
It suffices to note that
the limit of a sequence in ${\bf T}$ is
path-connected (see, for example, Theorem~7.5.1 in \cite{BurBurIva01}) and 
satisfies the four point
condition (\ref{four}), (indeed, as remarked after Proposition~7.4.12 
in  \cite{BurBurIva01}, there is a ``meta--theorem'' that
if a feature of a compact metric space
can be formulated as a continuous property of
distances among finitely many points, then this feature 
is preserved under Gromov-Hausdorff limits).
\end{proof}

\bigskip

\noi
{\em Proof of Theorem~\ref{T0} }\quad 
We start by showing {\em separability}.
Given a compact $\R$-tree, $T$, and $\varepsilon>0$, let
$S_{\varepsilon}$ be a finite $\varepsilon$-net in $T$. 
For $a,b\in T$, let 
\be{arc}
   [a,b[\,:=\phi_{a,b}(\,[0,d(a,b)[\,)\quad \mbox{and}\quad  
   ]a,b[\,:=\phi_{a,b}(\,]0,d(a,b)[\,) 
\ee
be the unique half open and open, 
respectively, {\em arc} between them, and
write $T_\varepsilon$ for the {\em subtree of $T$ spanned by $S_\varepsilon$},
that is,
\be{Teps}
   T_\varepsilon
 :=
   \bigcup\nolimits_{x,y\in S_{\varepsilon}}[x,y] \quad 
   \mbox{and}\quad  d_{T_\varepsilon}:=d\big|_{T_\varepsilon}.
\ee
Obviously, $T_\varepsilon$ is still an $\varepsilon$-net for $T$,
and hence
$d_{{\mathrm{GH}}}(T_\varepsilon,T)\le d_{H}(T_\varepsilon,T)\le\varepsilon$.

Now each $T_\varepsilon$ is just a ``finite tree with edge-lengths''
and can clearly be approximated arbitrarily closely in the 
$d_{{\mathrm{GH}}}$-metric 
by trees with the same tree topology (that is, ``shape''), and 
rational edge-lengths.
The set of isometry types of finite trees with rational edge-lengths
is countable, and so $({\bf T},d_{{\mathrm{GH}}})$ is separable.\sm

It remains to establish {\em completeness}.  
It suffices by Lemma~\ref{L:Tclosed} to show that
any Cauchy sequence in ${\bf T}$ converges
to some compact metric space, or, equivalently, any
Cauchy sequence in ${\bf T}$ has a subsequence that converges
to some metric space.

Let $(T_n)_{n\in\N}$ be a Cauchy
sequence in ${\bf T}$. By Exercise~7.4.14 and
Theorem~7.4.15 in \cite{BurBurIva01},
a sufficient condition for this sequence to have a subsequential
limit is that for every $\varepsilon>0$ there exists a
positive number $N=N(\varepsilon)$ such that every $T_n$ contains an 
$\varepsilon$-net of cardinality $N$.


Fix $\varepsilon>0$ and $n_0=n_0(\varepsilon)$ such that
$d_{{\mathrm{GH}}}(T_m,T_n)<\varepsilon/2$ for $m,n\ge n_0$. Let $S_{n_0}$ be a finite
$(\varepsilon/2)$-net for $T_{n_0}$ of cardinality $N$. 
Then by
(\ref{GH}) for each $n\ge n_0$ there exists a correspondence ${\Re}_n$ 
between $T_{n_0}$ and $T_n$ such that
${\mathrm{dis}}(\Re_n)<\varepsilon$. For each $x\in
T_{n_0}$, choose $f_n(x)\in T_n$ such that $(x,f_n(x))\in{\Re_n}$. 
Since for any $y\in T_n$ with $(x,y)\in{\Re_n}$,
$d_{T_n}(y,f_n(x))\le{\mathrm{dis}}(\Re_n)$, 
for all $n\ge n_0$, the set $f_n(S_{n_0})$
is an $\varepsilon$-net of cardinality $N$
for $T_{n}$, $n\ge n_0$. \hfill $\qed$ \bi

\subsection{Unrooted trees with $4$ leaves}
\label{fourp}
For the sake of reference and establishing some notation,
we record here some well-known facts about reconstructing trees from 
a knowledge of the distances between the leaves.  We remark that
the fact that trees can be reconstructed from their collection of
leaf-to-leaf distances (plus also the leaf-to-root distances
for rooted trees) is of huge practical importance in so-called
{\em distance methods} for inferring phylogenetic trees from DNA
sequence data, and the added fact that one can build such trees
by building subtrees for each collection of four leaves
is the starting point for the sub-class of
distance methods called {\em quartet methods}. We refer the reader
to \cite{Fel03, SemStee03} for an extensive description of these
techniques and their underlying theory.

\begin{lemma}\label{4fact} The isometry class
of an unrooted tree $(T,d)$ with four leaves
is uniquely determined by the distances between the leaves of $T$.
\end{lemma}

\begin{proof} Let $\{x_1,x_2,x_3,x_4\}$ be the set of 
leaves of $T$. The tree $T$ has one of four possible shapes: 

\setlength{\unitlength}{0.7pt}
\begin{picture}(300,170)(-20,-80)
\thicklines
\put(400,7){\line(-1,1){20}}
\put(400,7){\line(1,1){30}}
\put(400,7){\line(-1,-1){21.1}}
\put(400,7){\line(1,-1){21.1}}
\put(360,-15){\mbox{$x_2$}}
\put(430,-15){\mbox{$x_4$}}
\put(360,25){\mbox{$x_1$}}
\put(430,45){\mbox{$x_3$}}
\put(350,55){\mbox{(IV)}}
\put(280,32){\line(-1,1){20}}
\put(280,32){\line(1,1){30}}
\put(280,32){\line(0,-1){25}}
\put(280,7){\line(-1,-1){21.1}}
\put(280,7){\line(1,-1){21.1}}
\put(240,-15){\mbox{$x_2$}}
\put(310,-15){\mbox{$x_3$}}
\put(250,35){\mbox{$x_1$}}
\put(310,45){\mbox{$x_4$}}
\put(220,55){\mbox{(III)}}
\put(150,32){\line(-1,1){20}}
\put(150,32){\line(1,1){30}}
\put(150,32){\line(0,-1){25}}
\put(150,7){\line(-1,-1){21.1}}
\put(150,7){\line(1,-1){21.1}}
\put(110,-15){\mbox{$x_2$}}
\put(180,-15){\mbox{$x_4$}}
\put(120,35){\mbox{$x_1$}}
\put(180,45){\mbox{$x_3$}}
\put(90,55){\mbox{(II)}}
\put(20,32){\line(-1,1){20}}
\put(20,32){\line(1,1){30}}
\put(20,32){\line(0,-1){25}}
\put(20,7){\line(-1,-1){21.1}}
\put(20,7){\line(1,-1){21.1}}
\put(-20,-15){\mbox{$x_1$}}
\put(50,-15){\mbox{$x_2$}}
\put(-10,35){\mbox{$x_4$}}
\put(55,55){\mbox{$x_3$}}
\put(30,10){\mbox{$y_{1,2}$}}
\put(30,30){\mbox{$y_{3,4}$}}
\put(-30,55){\mbox{(I)}}
\put(-20,-60){\makebox{{\sc Figure~2} 
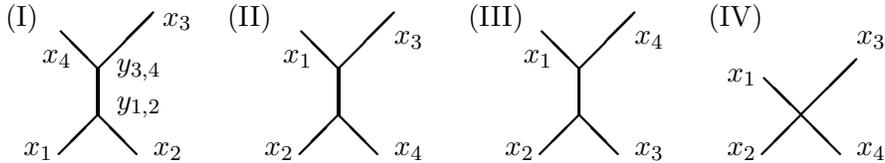
{\em shows the $4$ different
shapes of a labeled tree with $4$ leaves.}}}
\end{picture}

Consider case
$(I)$, and let $y_{1,2}$ be the uniquely determined branch
point on the tree that lies on the arcs $[x_1,x_2]$ and $[x_1,x_3]$,
and $y_{3,4}$ be the uniquely determined branch
point on the tree that lies on the arcs $[x_3,x_4]$ and $[x_1,x_3]$. 
Observe that
\be{conc}
\begin{aligned}
   d(x_1,y_{1,2})&=\frac{1}{2}(d(x_1,x_2)+d(x_1,x_3)-d(x_2,x_3))
  \\
   d(x_2,y_{1,2})&=\frac{1}{2}(d(x_1,x_2)+d(x_2,x_3)-d(x_1,x_3))
  \\
   d(x_3,y_{3,4})&=\frac{1}{2}(d(x_3,x_4)+d(x_1,x_3)-d(x_1,x_4))
  \\
   d(x_4,y_{3,4})&=\frac{1}{2}(d(x_3,x_4)+d(x_1,x_4)-d(x_1,x_3))
  \\
   d(y_{1,2},y_{3,4})&=\frac{1}{2}(d(x_1,x_4)+d(x_2,x_3)-
   d(x_1,x_2)-d(x_3,x_4)),
\end{aligned}
\ee   

Similar observations for the other cases show that if we know the
shape of the tree, then we can determine its edge-lengths from
leaf-to-leaf distances. Note also that
\be{char}
\begin{aligned}
   \chi_{(I)}(T)
 &:=
   \frac{1}{2}\left(d(x_1,x_3)+d(x_2,x_4)-d(x_1,x_2)-d(x_3,x_4)\right)
  \\[1mm]
 &\left\{\begin{array}{lc}>0 & \mbox{ \quad for shape (I)}, \\[1mm]
 <0 & \mbox{ \quad for shape (II)}, \\[1mm]
 =0 & \mbox{ \quad for shapes (III) and (IV)}\end{array}\right..
\end{aligned}
\ee\sm

This and analogous inequalities for the quantities that reconstruct
the length of the ``internal'' edge in shapes $(II)$ and $(III)$, 
respectively, show that the shape of the tree can also be
reconstructed from leaf-to-leaf distances. 
\end{proof}

\subsection{Rooted $\R$-trees}
\label{roottree}
Since we are mainly interested in rooted trees, we extend our 
definition as follows: A {\em rooted $\R$-tree}, $(X,d,\rho)$, is an
$\R$-tree $(X,d)$ with a distinguished point $\rho\in X$ that we
call the {\em root}.
It is helpful to use genealogical terminology
and think of $\rho$ as a common ancestor and $h(x):=d(\rho,x)$ 
as the real-valued generation to which 
$x\in X$ belongs ($h(x)$ is also called the {\em height} of $x$). 
We define a partial order 
$\le$ on $X$ by declaring (using the notation introduced
in (\ref{arc})) that $x\le y$ if 
$x\in [\rho,y]$, so that $x$ is an ancestor of $y$.
Each pair $x,y\in X$
has a well-defined {\em greatest common lower bound}, $x\wedge y$, in this 
partial order that we think of as the most recent common
ancestor of $x$ and $y$.

Let ${\bf T}^{\mathrm{root}}$ denote the collection of all root-invariant 
isometry classes
of rooted compact $\R$-trees, where we define a root-invariant 
isometry to be an isometry $\xi:         
(X_1,d_{X_1},\rho_1)\to(X_2,d_{X_2},\rho_2)$ with
$\xi(\rho_1)=\rho_2$.

We want to equip ${\bf T}^{\mathrm{root}}$ with a Gromov-Hausdorff
type distance that incorporates the special status of the
root.  We define the {\em rooted Gromov-Hausdorff distance},
$d_{{\mathrm{GH}}^{\mathrm{root}}}((X_1,\rho_1),(X_2,\rho_2))$, between
two rooted $\R$-trees $(X_1,\rho_1)$ and $(X_2,\rho_2)$ as the 
infimum of $d_{\mathrm H}(X_1^\prime,X_2^\prime)\vee
d_Z(\rho_1^\prime,\rho_2^\prime)$ 
over all rooted $\R$-trees
$(X_1^\prime,\rho_1^\prime)$ and $(X_2^\prime,\rho_2^\prime)$ 
that are root-invariant isomorphic to $(X_1,\rho_1)$ and $(X_2,\rho_2)$,
respectively,
and that are (as unrooted trees) subspaces of a common metric space 
$(Z,d_Z)$.\sm

As in (\ref{GH}), we can compute
 $d_{{\mathrm{GH}}^{\mathrm{root}}}((X_1,d_{X_1},\rho_1),(X_2,d_{X_2},\rho_2))$ 
by comparing distances  within $X_1$ to
distances within $X_2$, provided that the
distinguished status of the root is respected. 

\begin{lemma}
\label{L:GHroot}
For two 
rooted trees $(X_1,d_{X_1},\rho_1)$, and $(X_2,d_{X_2},\rho_2)$,
\be{GHroot}
   d_{{\mathrm{GH}}^{\mathrm{root}}}((X_1,d_{X_1},\rho_1),(X_2,d_{X_2},\rho_2))
 =
   \frac{1}{2}\inf_{{\Re^{\mathrm{root}}}}{\mathrm{dis}}({\Re^{\mathrm{root}}}),
\ee
where now the infimum is taken over all correspondences ${\Re^{\mathrm{root}}}$
between $X_1$ and $X_2$ with $(\rho_1,\rho_2)\in\Re^{\mathrm{root}}$.
\end{lemma}\bi

\begin{proof}
Indeed, for any root-invariant isometric copies
$(X^\prime_1,\rho^\prime_1)$ 
and $(X^\prime_2,\rho^\prime_2)$ embedded in $Z$, and 
$r>d_{{\mathrm{GH}}^{\mathrm{root}}}((X_1,\rho_1),(X_2,\rho_2))$,
\be{est14}
   {\Re}^{\mathrm{root}}
 :=
   \{(x_1,x_2);\,x_1\in X_1^\prime,x_2\in X_2^\prime,\,d_Z(x_1,x_2)<r\}
\ee
gives a correspondence 
between $X_1$ and $X_2$ containing 
$(\rho_1,\rho_2)$ such that ${\mathrm{dis}}(\Re^{\mathrm{root}})<2r$.  

On the other hand, given a correspondence
${\Re}^{\mathrm{root}}$
between $X_1$ and $X_2$ containing 
$(\rho_1,\rho_2)$, define a metric $d_{X_1\amalg X_2}$ on 
the disjoint union $X_1\amalg X_2$ 
such that the restriction of $d_{X_1\amalg X_2}$ to $X_i$
is $d_{X_i}$, for $i=1,2$, 
and for $x_1\in X_1$, $x_2\in X_2$, by
\be{est15}
\begin{aligned}
 & d_{X_1\amalg X_2}(x_1,x_2)
  \\
 & \quad :=
   \inf\{d_{X_1}(x_1,y_1)+d_{X_2}(x_2,y_2)+
      \frac{1}{2} {\mathrm{dis}}(\Re^{\mathrm{root}}) \, :\,
   (y_1,y_2)\in{\Re}^{\mathrm{root}}\}
\end{aligned}
\ee 
-- in particular, if the pair
$(x_1,x_2)$ actually belongs to the
correspondence $\Re^{\mathrm{root}}$, then $d_{X_1\amalg
  X_2}(x_1,x_2)=\frac12dis(\Re^{\mathrm{root}})$.
We leave it to the reader to check that $d_{X_1\amalg X_2}$
is, indeed, a metric.
Then (computing Hausdorff distance within  $X_1\amalg X_2$
using  $d_{X_1\amalg X_2}$)
we have 
\be{est15a}
d_{\mathrm H}(X_1,X_2)\vee d_{X_1\amalg X_2}(\rho_1,\rho_2)\le
\frac{1}{2} {\mathrm{dis}}(\Re^{\mathrm{root}}). 
\ee
\end{proof}

We state an analogue of Theorem~\ref{T0} for rooted compact $\R$-trees.

\begin{theorem}\label{T3} The metric space
$({\bf T}^{\mathrm{root}},d_{{GH}^{\mathrm{root}}})$ is 
  complete and separable.
\end{theorem}\sm

Before we can prove Theorem~\ref{T3}
we need two preparatory
results. The first is a
the counterpart of Corollary~7.3.28 in \cite{BurBurIva01} and
presents convenient upper and
lower estimates for $d_{{\mathrm{GH}}^{\mathrm{root}}}$
that differ by a multiplicative constant.

Let $(X_1,\rho_1)$ and $(X_2,\rho_2)$ be two rooted compact $\R$-trees,
and take $\varepsilon>0$. 
A map $f$ is called a {\em root-invariant
$\varepsilon$-isometry} from $(X_1,\rho_1)$ to $(X_2,\rho_2)$ if 
$f(\rho_1)=\rho_2$, 
${\mathrm{dis}}(f):= \sup\{|d_{X_1}(x,y) - d_{X_2}(f(x), f(y))| : x,y \in X_1\} <\varepsilon$ and $f(X_1)$
is an $\varepsilon$-net for $X_2$.

\begin{lemma}\label{7.3.28}
Let $(X_1,\rho_1)$ and $(X_2,\rho_2)$ be two rooted compact $\R$-trees,
and take $\varepsilon>0$. Then the following hold.
\begin{itemize}
\item[(i)] If $d_{{\mathrm{GH}}^{\mathrm{root}}}((X_1,\rho_1),(X_2,\rho_2))<\varepsilon$,
then there exists a root-invariant $2\varepsilon$-isometry from $(X_1,\rho_1)$
to $(X_2,\rho_2)$.
\item[(ii)] If there exists a root-invariant $\varepsilon$-isometry
  from $(X_1,\rho_1)$ to $(X_2,\rho_2)$, then 
$$
d_{{\mathrm{GH}}^{\mathrm{root}}}((X_1,\rho_1),(X_2,\rho_2))
\le\frac{3}{2}\varepsilon.
$$
\end{itemize}
\end{lemma}

\begin{proof}
(i) Let $d_{{\mathrm{GH}}^{\mathrm{root}}}((X_1,\rho_1),(X_2,\rho_2))<\varepsilon$. By 
 Lemma~\ref{L:GHroot} there exists a correspondence $\Re^{\mathrm{root}}$
between $X_1$ and $X_2$ such that $(\rho_1,\rho_2)\in\Re^{\mathrm{root}}$
and ${\mathrm{dis}}(\Re^{\mathrm{root}})<2\varepsilon$. Define $f:X_1\to X_2$
by setting $f(\rho_1)=\rho_2$, and 
choosing $f(x)$ such that $(x,f(x))\in\Re^{\mathrm{root}}$ for all 
$x\in X_1\setminus\{\rho_1\}$. Clearly, 
${\mathrm{dis}}(f)\le{\mathrm{dis}}(\Re^{\mathrm{root}})<2\varepsilon$. To see that
$f(X_1)$ is an $2\varepsilon$-net for $X_2$, let $x_2\in X_2$,
and choose $x_1\in X_1$ such that $(x_1,x_2)\in\Re^{\mathrm{root}}$.
Then $d_{X_2}(f(x_1),x_2)\le
d_{X_1}(x_1,x_1)+{\mathrm{dis}}(\Re^{\mathrm{root}})<2\varepsilon$.\sm

(ii) Let $f$ be a root-invariant $\varepsilon$-isometry
  from $(X_1,\rho_1)$ to $(X_2,\rho_2)$. Define a correspondence
  $\Re^{\mathrm{root}}_f\subseteq X_1\times X_2$ by
\be{ccc}
   \Re^{\mathrm{root}}_f:=\{(x_1,x_2):\,d_{X_2}(x_2,f(x_1))\le\varepsilon\}.
\ee 
Then $(\rho_1,\rho_2)\in\Re^{\mathrm{root}}_f$ and 
$\Re^{\mathrm{root}}_f$ is indeed a correspondence since $f(X_1)$ is a
  $\varepsilon$-net for $X_2$. If
  $(x_1,x_2),(y_1,y_2)\in\Re^{\mathrm{root}}_f$,
then
\be{diis}
\begin{split}
   |d_{X_1}(x_1,y_1)-d_{X_2}(x_2,y_2)|
 &\le 
   |d_{X_2}(f(x_1),f(y_1))-d_{X_1}(x_1,y_1)| \\
 & \quad  +d_{X_2}(x_2,f(x_1))+d_{X_2}(f(x_1),y_2) \\
 & <3\varepsilon. \\
\end{split}
\ee
Hence ${\mathrm{dis}}(\Re^{\mathrm{root}}_f)<3\varepsilon$ and, by (\ref{GHroot}),
$d_{{\mathrm{GH}}^{\mathrm{root}}}((X_1,\rho_1),(X_2,\rho_2))\le\frac{3}{2}\varepsilon$.
\end{proof}

The second preparatory result we need is the following compactness criterion, which
is the analogue of Theorem~7.4.15 in \cite{BurBurIva01} (note also
Exercise~7.4.14 in \cite{BurBurIva01}) and can be proved the same
way, using Lemma~\ref{7.3.28} in place of Corollary~7.3.28 in \cite{BurBurIva01}
and noting that the analogue of Lemma~\ref{L:Tclosed} holds for ${\bf T}^{\mathrm{root}}$.

\begin{lemma}
\label{L:pointcompactcrit}
A subset ${\mathcal T} \subset {\bf T}^{\mathrm{root}}$ 
is pre-compact if for every 
$\varepsilon>0$
there exists a positive integer $N(\varepsilon)$ such that 
each $T \in {\mathcal T}$ has an
$\varepsilon$-net with at most
$N(\varepsilon)$ points.
\end{lemma}

\noi  
{\em Proof of Theorem~\ref{T3} }\quad 
The proof follows very much the same lines as that of Theorem~\ref{T0}.
The proof of {\em separability} is almost identical.
The key step in establishing {\em completeness} is
again to show that a Cauchy
sequence in ${\bf T}^{\mathrm{root}}$ has
a subsequential limit.  This can be shown in the same manner
as in the proof of Theorem~\ref{T0}, with 
an appeal to Lemma~\ref{L:pointcompactcrit} 
replacing  one to Theorem~7.4.15 and
Exercise~7.4.14 in \cite{BurBurIva01}.
\hfill $\qed$ \bi\bi

\subsection{Length measure}

Recall that the root growth with re-grafting dynamics
involve points being chosen uniformly at random on
a finite tree.  In order to extend the dynamics to
general rooted compact $\R$-trees, we will
require the fact that rooted
compact $\R$-trees are associated with 
a natural length measure as follows. 
Fix $(T,d,\rho)\in{\bf T}^{\mathrm{root}}$,
and denote the Borel-$\sigma$-field on $T$ by ${\mathcal B}(T)$.
For $a,b\in T$, recall the half open arc $[a,b[$ from (\ref{arc}), and
let
\be{sce}
   T^o := \bigcup_{b\in T}\,[\rho,b[
\ee
the {\em skeleton} of $T$. Observe that if 
$T^\prime \subset T$ is a
dense countable set, then (\ref{sce}) holds with $T$ replaced by
$T^\prime$. In particular, $T^o \in {\mathcal B}(T)$ and
${\mathcal B}(T)\big|_{T^o}=\sigma(\{]a,b[;\,a,b\in T^\prime\})$, where
${\mathcal B}(T)\big|_{T^o}:=\{A \cap T^o;\,A\in{\mathcal B}(T)\}$. Hence there
exist a unique $\sigma$-finite measure $\mu$ on $T$, called {\em length
  measure}, such that $\mu(T\setminus T^o)=0$ and 
\be{length}
   \mu(]a,b[)=d(a,b),\quad \forall\;a,b\in T.
\ee
In particular, $\mu$ is the trace onto $T^o$ of one-dimensional 
Hausdorff measure on $T$. 

\bigskip
\noi
{\em Remark} \quad  The terminology {\em skeleton} might seem
somewhat derisory, since for finite trees the difference between
the skeleton and the whole tree is just a finite number of points.
However, it is not difficult to produce $\R$-trees for which
the difference between the skeleton and the whole tree is
a set with Hausdorff dimension greater than one (the Brownian
CRT will almost surely be such a tree).  This explains
our requirement that $\mu$ is carried by the skeleton.

\bigskip
\noi
{\em Remark} \quad  Elements of ${\bf T}^{\mathrm{root}}$
are really equivalence classes of trees rather than trees themselves,
so what we are describing here is a way of associating a measure
to each element of the equivalence class.  However, this procedure
respects the equivalence relation in that if $T'$ and $T''$
are two representatives of the same equivalence class and
are related by a root-invariant isometry $f:T' \rightarrow T''$,
then the associated length measures $\mu'$ and $\mu''$
are such that $\mu''$ is the push-forward of $\mu'$ by $f$
and $\mu'$ is the push-forward of $\mu''$ by the inverse of $f$
(that is, $\mu''(A'') = \mu'(f^{-1}(A''))$ and 
$\mu'(A') = \mu''(f(A'))$ for Borel sets $A'$ and $A''$
of $T'$ and $T''$, respectively).

\subsection{Rooted subtrees and trimming}

Recall from the Introduction that our strategy for extending
the root growth with re-grafting dynamics from finite trees
will involve a limiting procedure in which a general
rooted compact $\R$-tree is approximated ``from the inside''
by an increasing sequence of finite subtrees.  We therefore need to
establish some facts about such approximations.

To begin with, we require a notation for one tree being a subtree
of another, with both trees sharing the same root.  We need to
incorporate the fact that we are dealing with
equivalence classes of trees rather than trees
themselves.  A {\em rooted subtree} of 
$(T,d,\rho)
\in{\bf T}^{\mathrm{root}}$
is an element
$(T^{\ast},d^\ast,\rho^{\ast}),
\in{\bf T}^{\mathrm{root}}$
that has a class representative that is
a subspace of a class representative of
$(T,d,\rho)$, with the
two roots coincident.  Equivalently,
any class representative of $(T^{\ast},d^\ast,\rho^{\ast})$ 
can be isometrically embedded into any class
representative of $(T,d,\rho)$ via an isometry that maps
roots to roots.  We write $T^{\ast}\preceq^{\mathrm{root}} T$
and note that $\preceq^{\mathrm{root}}$ is an partial
order on ${\bf T}^{\mathrm{root}}$.



All of the ``wildness'' in a compact $\R$-tree happens ``at the leaves''.
For example, if $T \in {\bf T}^{\mathrm{root}}$ has a point
$x$ at which infinite branching occurs (so that the removal of
$x$ would disconnect $T$ into infinitely many components),
then any open neighborhood of $x$ must contain infinitely many
leaves, while for each $\eta>0$ there are only finitely many leaves $y$ such
that $x \in [\rho,y]$ with $d(x,y)>\eta$.  
A natural way in which to produce a finite subtree
that approximates a given tree is thus to fix $\eta>0$
and trim off the fringe
of the tree by removing those points that are not at least
distance $\eta$ from at least one leaf.
Formally, for $\eta>0$ define 
$R_\eta:{\bf T}^{\mathrm{root}}\to{\bf T}^{\mathrm{root}}$ to be 
the map that assigns to $(T,\rho)\in{\bf T}^{\mathrm{root}}$ the 
rooted subtree
$(R_\eta(T),\rho)$ 
that consists 
of $\rho$ and points $a \in T$ for which the subtree
\be{subab}
   S^{T,a}:=\{x\in T:\,a\in[\rho,x[\} 
\ee
(that is, the {\em subtree above $a$})
has height greater than
or equal to $\eta$. Equivalently,
\be{trim}
   R_\eta(T)
 :=
   \{x\in T:\,\exists\,y\in T\;x\in[\rho,y],\,d_T(x,y)\ge\eta\} \cup \{\rho\}.
\ee
In particular, if $T$ has height at most $\eta$, then
$R_\eta(T)$ is just the trivial tree consisting of the root $\rho$.
\sm
 
\noi{\em Remark }\quad 
Notice that the map described in (\ref{trim}) maps a metric space into
a sub-space. However, since isometric spaces are mapped into isometric
sub-spaces, we may think of $R_\eta$ as a map from ${\bf T}^{\mathrm{root}}$
into ${\bf T}^{\mathrm{root}}$.\bi  

\begin{lemma}\label{L1}
\begin{itemize}
\item[(i)]
The range of $R_{\eta}$
consists of finite rooted trees.
\item[(ii)]
The map $R_{\eta}$ is continuous.
\item[(iii)]
The family of maps $(R_{\eta})_{\eta>0}$
is a semigroup; that is, 
$R_{\eta'} \circ R_{\eta''} = R_{\eta'+\eta''}$ for
$\eta',\eta''>0$.  In particular, 
$R_{\eta'}(T) \preceq^{\mathrm{root}}  R_{\eta''}(T)$
for $\eta' \ge \eta'' > 0$.
\item[(iv)] 
For any $(T,\rho) \in {\bf T}^{\mathrm{root}}$, 
$d_{{\mathrm{GH}}^{\mathrm{root}}}((T,\rho),(R_\eta(T),\rho)) 
\le d_{\mathrm{H}}(T, R_\eta(T)) \le \eta$.
\end{itemize}
\end{lemma}

\begin{proof}
(i) Fix $(T,d,\rho)\in{\bf T}^{\mathrm{root}}$. Let 
$E\subset R_{\eta}(T)$ be the leaves of 
$R_{\eta}$, that is, the points that have no subtree above them. We
have to show that $E$ is finite. However, if $a_1,a_2,\ldots$ are
infinitely many points in $E\setminus\{\rho\}$,   
then we can find points $b_1,b_2,\ldots$ in $T$ such that $b_i$ is in the
subtree above $a_i$ and $d(a_i,b_i)\ge\eta$. It follows that
$\inf_{i\not =j}d(b_i,b_j)\ge 2\eta$, which contradicts the compactness of
$T$.\sm

\noi
(ii) Suppose that $(T^\prime,d^\prime,\rho^\prime)$  
and $(T^{\prime\prime},d^{\prime\prime},\rho^{\prime\prime})$ are two
compact trees with
\begin{equation*}
d_{{\mathrm{GH}}^{\mathrm{root}}}((T^\prime,\rho^\prime),(T^{\prime\prime},\rho^{\prime\prime}))<
\varepsilon.
\end{equation*} 
By part~(i) of  Lemma~\ref{7.3.28}
there exists a root-invariant $2\epsilon$-isometry $f:T^\prime\to T^{\prime\prime}$.
Recall that this means, $f(\rho^\prime) = \rho^{\prime\prime}$,
${\mathrm{dis}}(f)<2\varepsilon$, and
$f(T^\prime)$ is an $2\varepsilon$-net for $T^{\prime\prime}$.

For $a\in R_{\eta}(T^\prime)$, let $\bar{f}(a)$ be the unique point in 
$R_{\eta}(T^{\prime\prime})$ that is closest to $f(a)$. 
We will show that
$\bar f: R_{\eta}(T^\prime) \rightarrow R_{\eta}(T^{\prime\prime})$
is a root-invariant $25 \varepsilon$-isometry and hence, by part~(ii) of Lemma~\ref{7.3.28},
$d_{{\mathrm{GH}}^{\mathrm{root}}}(R_{\eta}(T^\prime), R_{\eta}(T^{\prime\prime})) \le \frac{3}{2} 25 \varepsilon$.

We first show that
\be{est3}
   \sup\{d^{\prime\prime}(f(a),\bar{f}(a)) : a\in R_{\eta}(T^\prime)\}
 \le 
   8\varepsilon.
\ee 
Fix $a\in R_{\eta}(T^\prime)$
and let $b\in T^{\prime}$ be a point in the subtree above $a$ such
that $d^\prime(a,b)\ge\eta$. Denote the most recent common ancestor 
of $f(a)$ and $f(b)$ on $T^{\prime\prime}$ by 
$f(a)\wedge^{\prime\prime}f(b)$. 
 
\setlength{\unitlength}{0.7pt}
\begin{picture}(300,190)(-20,-100)
\thinlines
\put(220,-28){\line(0,1){90}}
\put(220,-28){\circle{3}}
\put(20,-28){\line(0,1){80}}
\put(220,18){\line(1,0){20}}
\put(20,-28){\circle{3}}
\put(20,52){\circle{3}}
\put(20,23){\circle{3}}
\put(220,18){\circle{3}}
\put(220,62){\circle{3}}
\put(240,18){\circle{3}}
\put(135,17){\mbox{$f(a)\wedge^{\prime\prime} f(b)$}}
\put(230,-25){\mbox{$\rho^{\prime\prime}$}}
\put(242,15){\mbox{$f(a)$}}
\put(230,57){\mbox{$f(b)$}}
\put(30,-25){\mbox{$\rho^{\prime}$}}
\put(30,17){\mbox{$a$}}
\put(30,47){\mbox{$b$}}
\put(-20,-60){\makebox{
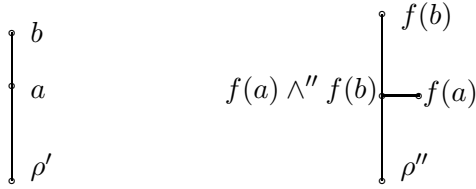
\begin{minipage}[t]{12cm}{\em {\sc Figure~3} illustrates the shapes
of the trees spanned by $\{\rho^\prime,a,b\}$ and by
$\{\rho^{\prime\prime},f(a),f(b)\}$. The point $\bar{f}(a)$
lies somewhere on the arc $[\rho^{\prime\prime},f(a)]$.}\end{minipage}}}
\end{picture}

Then
\be{est4}
\begin{aligned}
   d^{\prime\prime}&\left(f(a)\wedge^{\prime\prime}f(b),f(a)\right)
  \\[1mm]
 &=
   \frac{1}{2}\left(d^{\prime\prime}(f(a),f(b))+
   d^{\prime\prime}(\rho^{\prime\prime},f(a))-
   d^{\prime\prime}(\rho^{\prime\prime},f(b))\right)
  \\[1mm]
 &\le
   \frac{1}{2}\left(\big|d^{\prime\prime}(f(a),f(b))-d^{\prime}(a,b)\big|\right.  \\[1mm]
 &\hspace{1.5cm}\left.+
   \big|d^{\prime\prime}(\rho^{\prime\prime},f(a))-d^{\prime}(\rho^{\prime},a)
   \big|+
   \big|d^{\prime\prime}(\rho^{\prime\prime},f(b))-
   d^{\prime}(\rho^{\prime},b)\big|\right)
  \\[1mm]
 &\le
   3\varepsilon.
\end{aligned}
\ee 

If $\bar{f}(a)\in\,[f(a)\wedge^{\prime\prime}f(b),f(a)]$ 
then we are immediately done. 
Otherwise, $\bar{f}(a)\in\,[\rho^{\prime\prime},f(a)]$ and
$\bar{f}(a)$ is a leaf in
$R_{\eta}(T^{\prime\prime})$. Hence $f(b)\not\in
R_\eta(T^{\prime\prime})$, 
and therefore
\be{est4a}
   d^{\prime\prime}(\bar{f}(a),f(b)) \le \eta.
\ee
Furthermore,
\be{est5}
\begin{aligned}
   d^{\prime\prime}\left(f(a)\wedge^{\prime\prime}\hspace{-.1cm}f(b),
   f(b)\right)
 &=
   d^{\prime\prime}(f(a),f(b))-d^{\prime\prime}
   \left(f(a)\wedge^{\prime\prime}\hspace{-.1cm}f(b),f(a)\right)
  \\[1mm]
 &\ge 
  \left(d^\prime(a,b)-2\varepsilon\right)-3\varepsilon
  \\[1mm]
 & \ge
  \eta-5\varepsilon.
\end{aligned}   
\ee 
Combining (\ref{est4}), (\ref{est4a}) and (\ref{est5}) finally yields that
$d^{\prime\prime}(\bar{f}(a),f(a)) \le 8\varepsilon$
and completes the proof of (\ref{est3}).

It follows from (\ref{est3}) that
\begin{equation*}
\begin{split}
{\mathrm{dis}}(\bar f) & = \sup\{|d'(a,b) - d''(\bar f(a), \bar f(b))| : a,b \in R_{\eta}(T^\prime)\} \\
& \le \sup\{|d'(a,b) - d''(f(a), f(b))| : a,b \in R_{\eta}(T^\prime)\} \\ 
& \quad + 2 \sup\{d^{\prime\prime}(f(a),\bar{f}(a)) : a\in R_{\eta}(T^\prime)\} \\
& < 2 \varepsilon + 2 \times 8 \varepsilon \\
& = 18 \varepsilon. \\
\end{split}
\end{equation*}

The proof of (ii) will thus be completed if we can show that $\bar f(R_{\eta}(T^\prime))$
is a $25 \varepsilon$-net in $R_{\eta}(T^{\prime\prime})$.
Consider a point $c \in R_{\eta}(T^{\prime\prime})$.  We need to show that there
is a point $b \in R_{\eta}(T^\prime)$ such that 
\be{netbarf}
d''(\bar f(b), c) < 25 \varepsilon.
\ee
If $d''(\rho'', c) < 7 \varepsilon$, then we are done, because we can take
$b = \rho'$ (recall that $\bar f(\rho') = \rho''$).  Assume, therefore, that
$d''(\rho'', c) \ge 7 \varepsilon$.  We can then find points $c_-, c_+ \in T^{\prime\prime}$
such that $\rho'' \le c_- \le c \le c+$ with $d''(c_-, c) = 7 \varepsilon$ and $d''(c,c_+) \ge \eta$.
There are corresponding points $a_-, a, a_+ \in T^\prime$ such that 
$d''(f(a_-),c_-) < 2 \varepsilon$,
$d''(f(a),c) < 2 \varepsilon$, and $d''(f(a_+),c_+) < 2 \varepsilon$.  We claim that
$b := a_- \wedge^\prime a_+$ (the most recent common ancestor
of $a_-$ and $a_+$ in the tree $T'$)
belongs to $R_{\eta}(T^\prime)$ and satisfies
(\ref{netbarf}).

Note first of all that
\begin{equation*}
\begin{split}
d'(b,a_+) & = d'(a_- \wedge^\prime a_+, a_+) \\
& = \frac{1}{2} \left(
d'(a_+, a_-) + d'(\rho', a_+) - d'(\rho', a_-) \right) \\
& \ge \frac{1}{2} \bigl(
d''(f(a_+), f(a_-)) -2 \varepsilon + d''(f(\rho'), f(a_+)) -2 \varepsilon \\
& \quad - d''(f(\rho'), f(a_-)) - 2 \varepsilon \bigr) \\
& \ge \frac{1}{2} \left(
d''(c_+, c_-) -4 \varepsilon + d''(\rho'', c_+) -2 \varepsilon - d''(\rho'', c_-) - 2 \varepsilon \right) - 3 \varepsilon \\
& = d''(c_+, c_-) - 7 \varepsilon \\
& = \eta + 7 \varepsilon - 7 \varepsilon \\
& \eta, \\
\end{split}
\end{equation*}
and so $b \in R_{\eta}(T^\prime)$.

Furthermore,
\begin{equation*}
\begin{split}
d''(c, f(b)) & \le d''(c,c_-) + d''(c_-, f(a_-)) + d''(f(a_-), f(b)) \\
& \le 7 \varepsilon + 2 \varepsilon + d'(a_-, b) + 2 \varepsilon \\
& = 11 \varepsilon + \frac{1}{2} \left(
d'(a_+, a_-) + d'(\rho', a_-) - d'(\rho', a_+)\right) \\
& \le 11 \varepsilon + \frac{1}{2} \bigl(
d''(f(a_+), f(a_-)) + 2 \varepsilon + d''(f(\rho'), f(a_-)) + 2 \varepsilon  \\
& \quad - d'(f(\rho'), f(a_+)) + 2 \varepsilon \bigr) \\
& \le 14 \varepsilon + \frac{1}{2} \left(
d''(c_+, c_-) + 2 \varepsilon + d''(\rho'', c_-) +2 \varepsilon - d''(\rho'', c_+) +2 \varepsilon \right) \\
& = 17 \varepsilon. \\
\end{split}
\end{equation*}
Therefore, by (\ref{est3}),
\begin{equation*}
d(c, \bar f(b)) \le 17 \varepsilon + 8 \varepsilon = 25 \varepsilon.
\end{equation*}
This completes the proof of (\ref{netbarf}), and thus the proof of part~(ii).

\noi
Claims (iii) and (iv) are clear.
\end{proof}


Finally, we require the following result, which will be the key
to showing that the ``projective limit'' of a consistent
family of tree-valued processes can actually be thought
of as a tree-valued process in its own right. 

\begin{lemma}\label{L00} 
Consider a sequence 
$(T_n)_{n \in \N}$ of representatives of isometry classes
of rooted compact trees in $({\bf T},d_{{\mathrm{GH}}^{\mathrm{root}}})$ 
with the following properties.
\begin{itemize}
\item
Each set $T_n$ is a subset of some common set $U$.
\item
Each tree $T_n$ has the same root $\rho \in U$.
\item
The sequence $(T_n)_{n \in \N}$ is nondecreasing, that is, 
$T_1\subseteq T_2\subseteq\cdots\subseteq U$.
\item
Writing $d_n$ for the metric on $T_n$, for $m<n$ 
the restriction of $d_n$ to  $T_m$ coincides
with $d_m$, so that there is a well-defined metric on
$T:=\bigcup_{n\in\N} T_n$ 
given by
\be{me}
   d(a,b)=d_n(a,b),\quad a,b\in T_n.
\ee 
\item The sequence of subsets $(T_n)_{n\in\N}$ 
is Cauchy in the Hausdorff
distance with respect to $d$.
\end{itemize}
Then the following hold.
\begin{itemize}
\item[(i)]
The metric completion 
$\bar{T}$ of $T$  is a compact $\R$-tree,
and $d_{\mathrm H}(T_n,\bar{T}) \rightarrow 0$
as $n \rightarrow \infty$, where the Hausdorff 
distance is computed with
respect to the extension of $d$ to $\bar{T}$. In particular,
\be{mee}
\lim_{n\to\infty}d_{{\mathrm{GH}}^{\mathrm{root}}}((T_n,\rho),(\bar{T},\rho))=0.
\ee
\item[(ii)]
The tree $\bar{T}$ has skeleton $\bar{T}^o = \bigcup_{n \in \N} T_n^o$.
\item[(iii)]  The length measure on $\bar{T}$ is the unique measure
concentrated on $\bigcup_{n \in \N} T_n^o$ that restricts to the length
measure on $T_n$ for each $n \in \N$.
\end{itemize}
\end{lemma}

\begin{proof}
(i) Because $\bar{T}$ is a complete metric space, 
the collection of closed subsets
of $\bar{T}$ equipped with the Hausdorff distance is also complete
(see, for example, Proposition~7.3.7 of \cite{BurBurIva01}). 
Therefore the Cauchy sequence $(T_n)_{n \in \N}$ has a limit that is
(see Exercise~7.3.4 of \cite{BurBurIva01}) the closure of 
$\bigcup_{k\in\N}T_k$, i.e, $\bar{T}$ itself. 
It is clear that the complete space $\bar{T}$ is totally bounded, 
path-connected, and satisfies the four point condition, and
so $\bar{T}$ is a compact $\R$-tree. Finally,
\be{meep} 
   d_{{\mathrm{GH}}^{\mathrm{root}}}((T_n,\rho),(\bar{T},\rho))\le d_{\mathrm H}(T_n,\bar{T})\vee 
   d(\rho,\rho)=d_{\mathrm H}(T_n,\bar{T}) \rightarrow 0,
\ee
as $n \rightarrow \infty$.

Claims (ii) and (iii) are obvious.
\end{proof}

\section{Root growth with re-grafting}
\label{S:construct}

\subsection{Beginning the construction} 
\label{SS:beginconstruct}

We are now ready to begin in earnest the
construction of the ${\bf T}^{\mathrm{root}}$-valued 
Markov process, $X$, having the root growth with re-grafting
dynamics.   

Fix a tree $(T,d,\rho)\in{\bf T}^{\mathrm{root}}$. This
tree will be the initial state of $X$.  
We first recapitulate
the strategy outlined in the Introduction. 
In line with that semi-formal description, 
the ``stochastic inputs'' to the construction
of $X$ will be a collection of cut times and a corresponding
collection of cut points.

\begin{itemize}
\item
Construct simultaneously for each finite rooted subtree
$T^\ast\preceq^{\mathrm{root}} T$
a process $X^{T^\ast}$ with $X^{T^\ast}_0 = T^\ast$
that evolves according the root growth with re-grafting dynamics.
\item
Carry out this construction in such a way that if
$T^\ast$ and $T^{\ast \ast}$ are two finite subtrees with
$T^\ast\preceq^{\mathrm{root}} T^{\ast \ast}$, then 
$X^{T^\ast}_t \preceq^{\mathrm{root}} X^{T^{\ast \ast}}_t$
and the cut points for $X^{T^\ast}$ are those for $X^{T^{\ast \ast}}$
that happen to fall on $X^{T^\ast}_{\tau-}$ for a corresponding
cut time $\tau$ of  $X^{T^{\ast \ast}}$.  
Cut times $\tau$ for $X^{T^{\ast \ast}}$ 
for which the corresponding cut point does not fall on
$X^{T^\ast}_{\tau-}$ are not cut times for $X^{T^\ast}$.
\item 
The tree $(T,\rho)$ is a 
rooted Gromov-Hausdorff limit of finite $\R$-trees 
with root $\rho$ 
(indeed, any subtree spanned by a
finite $\varepsilon$-net and $\rho$ is finite
and has rooted Gromov-Hausdorff
distance less than $\varepsilon$ from $(T,\rho)$). 
In particular, $(T,\rho)$ is the ``smallest''
rooted compact $\R$-tree that contains all of the
finite rooted subtrees of $(T,\rho)$.
\item
Because of the consistent projective nature of the
construction, we can define $X_t:=X^{T}_t$ 
for $t \ge 0$ as the ``smallest''
element of
${\bf T}^{\mathrm{root}}$ that contains  
$X^{T^\ast}_t$, for all finite trees $T^\ast\preceq^{\mathrm{root}} T$.
\end{itemize}

It will be convenient for
establishing features of the process $X$ such as the strong Markov
property to introduce randomness later and work initially in a setting
where the cut times and cut points are fixed. 
There are two types of cut points: those that occur at
points which were present in the initial tree $T$
and those that occur at points which were added due to
subsequent root growth.  Accordingly, we
consider  two countable subsets $\pi_0 \subset \R^{++}\times T^o$ 
and $\pi \subset \{(t,x)\in\R^{++}\times\R^{++} : x \le t\}$.
(Once again we note that we are moving backwards 
and forwards between thinking of $T$
as a metric space or as an equivalence class of metric spaces.
As we have written things here, we are thinking of $\pi_0$
being associated with a particular class representative,
but of course $\pi_0$ corresponds to a similar set for any
representative of the same equivalence class
by mapping across using the appropriate
root invariant isometry.)

\begin{ass}
\label{orderliness}
Suppose that the sets $\pi_0$ and $\pi$ have the following properties.
\begin{itemize}
\item[(a)] For all $t_0>0$, each of the sets
$\pi_0 \cap (\{t_0\} \times T^o)$ and
$\pi \cap (\{t_0\} \times ]0,t_0])$ has at most one point
and at least one of these sets is empty.
\item[(b)] For all $t_0 >0$ and all finite subtrees
$T' \subseteq T$, the set
$\pi_0 \cap (]0,t_0] \times T')$ is finite.
\item[(c)] For all $t_0>0$, the set
$\pi \cap  \{(t,x)\in\R^{++}\times\R^{++} : x \le t \le t_0\}$
is finite.
\end{itemize}
\end{ass}

\bigskip

\noi
{\em Remark} \quad Conditions (a)--(c) of
Assumption \ref{orderliness}
will hold almost surely if $\pi_0$ and $\pi$
are realizations of Poisson point processes with respective
intensities $\lambda \otimes \mu$ and $\lambda \otimes \lambda$
(where $\lambda$ is Lebesgue measure),
and it is this random mechanism that we will introduce
later to produce a stochastic process having the
root growth with re-grafting dynamics.

\bigskip

It will be convenient to use the notations $\pi_0$ and
$\pi$ to also refer to the integer-valued measures that
are obtained by placing a unit point mass at each point
of the corresponding set.


Consider a finite rooted subtree $T^\ast\preceq^{\mathrm{root}} T$. 
It will avoid annoying circumlocutions about equivalence
via root invariant isometries if we work with particular
class representatives for $T^\ast$ and $T$, and, moreover,
suppose that $T^\ast$ is embedded in $T$.

Put $\tau^\ast_0:=0$, and let 
$0<\tau^\ast_1<\tau_2^\ast<\ldots$
(the cut times for $X^{T^\ast}$)
be the points of 
$\{t > 0:\,\pi_0(\{t\}\times T^\ast)>0\} 
\cup \{t > 0:\,\pi(\{t\}\times\R^{++})>0\}$.\sm

An explicit
construction of $X^{T^\ast}_t$ is then given 
in two steps:\bi

\noi
{\em Step~1 (Root growth) }\quad 
At any time $t\ge 0$, $X_t^{T^{\ast}}$ as a set is given by
the disjoint union $T^\ast\amalg ]0,t]$. The root of $X_t^{T^{\ast}}$
is the point $\rho_t:=t\in]0,t]$. The metric $d_t^{T^{\ast}}$ on $X_t^{T^{\ast}}$
is defined inductively as follows. Set $d_0^{T^\ast}$ to be the metric on 
$X_0^{T^\ast}=T^\ast$; that is, $d_0^{T^\ast}$ is the restriction of $d$
to $T^\ast$. Suppose that $d_t^{T^{\ast}}$ has been defined
for $0\le t\le\tau^\ast_n$. Define $d_t^{T^{\ast}}$ for $\tau_n^\ast<t<
\tau^\ast_{n+1}$ by
\be{rgr}
   d_t^{T^{\ast}}(a,b)
 :=
   \begin{cases}
   d_{\tau_n^\ast}(a,b), & \text{if  $a,b \in X^{T^\ast}_{\tau_n^\ast}$,}\\
   |b-a|, & \text{if $a,b\in]\tau_n^\ast,t]$,} \\
   |a-\tau_n^\ast|+d_{\tau_n^\ast}(\rho_{\tau_n^\ast},b), & 
   \text{if $a\in]\tau_n^\ast,t]$, $b\in X^{T^\ast}_{\tau_n^\ast}$.}
   \end{cases}
\ee

\noi
{\em Step~2 (Re-Grafting) }\quad
Note that the left-limit
$X^{T^\ast}_{\tau^\ast_{n+1}-}$ exists in the rooted Gromov-Hausdorff
metric.  As a set this left-limit is the disjoint union
\be{repr}
   X^{T^\ast}_{\tau^\ast_{n}}\amalg]\tau_n^\ast,\tau_{n+1}^\ast]\,
 =
   T^\ast\amalg ]0,\tau_{n+1}^\ast], 
\ee
and the corresponding metric $d_{\tau^\ast_{n+1}-}$ is
given by a prescription similar to (\ref{rgr}). 

\noi 
Define the $(n+1)^{\mathrm{st}}$ cut point for $X^{T^\ast}$ by
\be{cut}
   p^\ast_{n+1}
 :=
   \begin{cases}
     a\in T^\ast, & \text{if $\pi_0(\{(\tau^\ast_{n+1},a)\})>0$,}\\
     x\in]0,\tau^\ast_{n+1}], & \text{if $\pi(\{(\tau^\ast_{n+1},x)\})>0$.}
   \end{cases}
\ee

Let $S_{n+1}^\ast$ be the {\em subtree above $p_{n+1}^\ast$} in 
$X^{T^\ast}_{\tau^\ast_{n+1}-}$, that is,
\be{sub}
   S_{n+1}^\ast
 :=
   \{b\in X^{T^\ast}_{\tau^\ast_{n+1}-}:\,p_{n+1}^\ast\in\,[
   \rho_{\tau_{n+1}^\ast-},b[\,\}.
\ee

\noi 
Define the metric $d_{\tau_{n+1}^\ast}$ by
\be{pru}
\begin{split}
& d_{\tau_{n+1}^\ast}(a,b)\\
& :=
  \begin{cases}
     d_{\tau_{n+1}^\ast-}(a,b), & \text{if $a,b\in S_{n+1}^\ast$,}\\
     d_{\tau^\ast_{n+1}-}(a,b), & \text{if $a,b\in X^{T^\ast}_{\tau^\ast_{n+1}}
     \hspace{-.2cm}\setminus S_{n+1}^\ast$,}\\
     d_{\tau_{n+1}^\ast-}(a,\rho_{\tau^\ast_{n+1}})+
     d_{\tau_{n+1}^\ast-}(p_{n+1}^\ast,b), 
            & \text{if $a \in X^{T^\ast}_{\tau^\ast_{n+1}}
     \hspace{-.2cm}\setminus S_{n+1}^\ast$, $b \in  S_{n+1}^\ast$.} 
   \end{cases}\\
\end{split}
\ee

\noi
In other words $X^{T^\ast}_{\tau^\ast_{n+1}}$ is obtained from 
$X^{T^\ast}_{\tau^\ast_{n+1}-}$ by pruning off the subtree $S^\ast_{n+1}$
and re-attaching it to the root. \bi

Now consider two other finite, rooted subtrees $(T^{\ast\ast},\rho)$ and 
$(T^{\ast\ast\ast},\rho)$ of $T$ such that
$T^{\ast}\cup~T^{\ast\ast}\subseteq T^{\ast\ast\ast}$
(with induced metrics). Build
$X^{T^{\ast\ast}}$ and $X^{T^{\ast\ast\ast}}$ from $\pi_0$ and $\pi$
in the same
manner as $X^{T^\ast}$ (but starting at $T^{\ast\ast}$ and
$T^{\ast\ast\ast}$). It is clear from the construction that:
\begin{itemize}
\item
$X^{T^\ast}_t$
and $X_t^{T^{\ast\ast}}$ are rooted subtrees of
$X_t^{T^{\ast\ast\ast}}$ 
for all $t \ge 0$,
\item
the Hausdorff distance between $X_t^{T^{\ast}}$ and
$X_t^{T^{\ast\ast}}$ 
as subsets of $X_t^{T^{\ast\ast\ast}}$
does not depend on $T^{\ast\ast\ast}$, 
\item
the Hausdorff distance is constant between 
jumps of $X^{T^\ast}$ and $X^{T^{\ast\ast}}$ (when only root growth is occurring
in both processes). 
\end{itemize}
The following lemma shows that the Hausdorff distance between 
$X_t^{T^{\ast}}$ and $X_t^{T^{\ast\ast}}$ 
as subsets of $X_t^{T^{\ast\ast\ast}}$
does not increase at jump times.

\begin{lemma}\label{L0} 
Let $T$ be a finite rooted tree with root $\rho$ and metric $d$, and 
let $T^\prime$ and $T^{\prime\prime}$ be two 
rooted subtrees of $T$ 
(both with the induced metrics and root 
$\rho$). 
Fix $p\in T$, and let 
$S$ be the subtree in $T$
above $p$ (recall (\ref{sub})). 
Define a new metric $\hat{d}$ on $T$ by putting
\begin{equation*}
   \hat{d}(a,b)
 :=
   \begin{cases}
     d(a,b), & \text{if $a,b\in S$,}
  \\
     d(a,b), & \text{if $a,b\in T\setminus S$,}
  \\
     d(a,p)+d(\rho,b), & \text{if $a\in S,\,b\in T \setminus S$.}
   \end{cases}
\end{equation*}
Then the sets $T^\prime$ and $T^{\prime\prime}$ are also subtrees of
$T$ equipped with the induced metric $\hat{d}$, and the Hausdorff distance
between $T^\prime$ and $T^{\prime\prime}$ with respect to $\hat{d}$
is not greater than that with respect to $d$.
\end{lemma}\bi

%
%

\begin{proof}
Suppose that the Hausdorff distance
between $T^\prime$ and $T^{\prime\prime}$ under ${d}$ is less than
some given $\varepsilon>0$. Given $a\in T^\prime$, there then exists
$b \in T^{\prime\prime}$ such that $d(a,b)<\varepsilon$.  Because
$d(a,a \wedge b) \le d(a,b)$ and $a \wedge b \in T^{\prime\prime}$,
we may suppose (by replacing $b$ by $a \wedge b$ if necessary)
that $b \le a$.  We claim that $\hat d(a,c) < \varepsilon$
for some $c \in T^{\prime\prime}$.
This and the analogous result with the roles of 
$T^\prime$ and $T^{\prime\prime}$ interchanged will establish the result.

If $a,b\in S$
or $a,b\in T\setminus S$, then $\hat{d}(a,b)=d(a,b)<\varepsilon$. 
The only other possibility is that $a \in S$ and $b \in T \setminus S$,
in which case $p \in [b,a]$ (for $T$ equipped with $d$). Then
$\hat d(a,\rho) = d(a,p) \le d(a,b) < \varepsilon$, as required
(because $\rho \in T^{\prime\prime}$).
\end{proof}

Now let $T_1\subseteq T_2\subseteq \cdots$ be an increasing sequence 
of finite subtrees of $T$ such that $\bigcup_{n\in\N}T_n$ is dense in
$T$. Thus $\lim_{n\to\infty}d_{\mathrm H}(T_n,T)=0$. 
Let $X^1,X^2,\ldots$ be
constructed from $\pi_0$ and $\pi$ starting with
$T_1, T_2,\ldots$. Applying  Lemma~\ref{L0} yields
\be{pre}
   \lim_{m,n\to\infty}\sup_{t\ge 0}
       \,d_{{\mathrm{GH}}^{\mathrm{root}}}(X_t^m,X_t^n)=0.
\ee 
Hence by completeness of ${\bf T}^{\mathrm{root}}$, there exists a c{\`a}dl{\`a}g
${\bf T}^{\mathrm{root}}$-valued
process $X$ such that $X_0=T$ and
\be{cad}
   \lim_{m\to\infty}\sup_{t\ge 0} 
       \,d_{{\mathrm{GH}}^{\mathrm{root}}}(X_t^m,X_t)=0.
\ee  

{\em A priori,} the process $X$ could depend on the choice of the
approximating sequence of trees $(T_n)_{n \in \N}$.  To see that
this is not so, consider two approximating sequences
$T_1^1\subseteq T_2^1 \subseteq \cdots$ and
$T_1^2\subseteq T_2^2 \subseteq \cdots$.  For $k \in \N$,
write $T_n^3$ for the smallest rooted subtree of $T$ that contains
both $T_n^1$ and $T_n^2$.  As a set, $T_n^3 = T_n^1 \cup T_n^2$.
Now let $((X_t^{n,i})_{t \ge 0})_{n \in \N}$ for $i=1,2,3$
be the corresponding sequences of finite tree-value processes
and let $(X_t^{\infty,i})_{t \ge 0}$ for $i=1,2,3$ be the
corresponding limit processes.  By  Lemma~\ref{L0},
\be{approxnomatter}
\begin{split}
d_{ {\mathrm{GH}} ^{\mathrm{root}} }(X_t^{n,1}, X_t^{n,2})
& \le
d_{{\mathrm{GH}}^{\mathrm{root}}}(X_t^{n,1}, X_t^{n,3})
+
d_{{\mathrm{GH}}^{\mathrm{root}}}(X_t^{n,2}, X_t^{n,3}) \\
& \le
d_{\mathrm{H}}(X_t^{n,1}, X_t^{n,3})
+
d_{\mathrm{H}}(X_t^{n,2}, X_t^{n,3}) \\
& \le
d_{\mathrm{H}}(T_n^1, T_n^3)
+
d_{\mathrm{H}}(T_n^2, T_n^3) \\
& \le
d_{\mathrm{H}}(T_n^1, T)
+
d_{\mathrm{H}}(T_n^2, T) \rightarrow 0\\
\end{split}
\ee
as $n \rightarrow \infty$.  Thus, for each $t \ge 0$ the sequences
$(X_t^{n,1})_{n \in \N}$ and $(X_t^{n,2})_{n \in \N}$
do indeed have the same rooted Gromov-Hausdorff limit and the process
$X$ does not depend on the choice of approximating sequence
for the initial tree $T$.

\subsection{Finishing the construction}
\label{SS:finishconstruct}

In \ref{SS:beginconstruct} we constructed a 
${\bf T}^{\mathrm{root}}$-valued function $t \mapsto X_t$
starting with a fixed triple $(T, \pi_0, \pi)$, where 
$T \in {\bf T}^{\mathrm{root}}$ and $\pi_0, \pi$ satisfy the
conditions of Assumption \ref{orderliness}.  We now want
to think of $X$ as a function of time and such triples.

Let $\Omega^\ast$ be the set of triples $(T,\pi_0,\pi)$,
where $T$ is a rooted compact $\R$-tree (that is,
a class representative of an element of ${\bf T}^{\mathrm{root}}$)
and $\pi_0, \pi$ satisfy Assumption \ref{orderliness}.

The root invariant isometry equivalence relation on rooted compact
$\R$-trees extends naturally to an equivalence relation on
$\Omega^\ast$ by declaring that  two triples
$(T',\pi_0',\pi')$ and $(T'',\pi_0'',\pi'')$,
where $\pi_0' = \{(\sigma_i',x_i') : i \in \N\}$ and 
$\pi_0''= \{(\sigma_i'',x_i'') : i \in \N\}$,
are equivalent if there is a root invariant isometry $f$ mapping $T'$
to $T''$ and a permutation $\gamma$ of $\N$ 
such that $\sigma_i'' = \sigma_{\gamma(i)}'$ and
$x_i'' = f(x_{\gamma(i)}')$ for all $i \in \N$.
We write $\Omega$ for the resulting quotient space
of equivalence classes.

In order to do probability, we require that $\Omega$
has a suitable measurable structure.  We could do this
by specifying a metric on $\Omega$, but the following
approach is a little less cumbersome and suffices for our needs.  

Let
$\Omega^{\mathrm{fin}}$ denote the subset of $\Omega$
consisting of triples $(T,\pi_0,\pi)$ such that 
$T$, $\pi_0$ and $\pi$ are finite.  
We are going to define a metric on $\Omega^{\mathrm{fin}}$.
Let $(T',\pi_0',\pi')$ and $(T'',\pi_0'',\pi'')$
be two points in $\Omega^{\mathrm{fin}}$,
where $\pi_0'= \{(\sigma_1',x_1'),\ldots,(\sigma_p',x_p')\}$,
$\pi'=\{\tau_1', \ldots, \tau_r'\}$,
$\pi_0''= \{(\sigma_1'',x_1''),\ldots,(\sigma_q'',x_q'')\}$,
and
$\pi''=\{\tau_1'', \ldots, \tau_s''\}$.
Assume that
$0<\sigma_1' < \cdots < \sigma_p'$,
$0<\tau_1' < \cdots < \tau_r'$,
$0<\sigma_1''<\cdots<\sigma_q''$,
and
$0<\tau_1'' < \cdots < \tau_s''$.
The distance between $(T',\pi_0',\pi')$ and $(T'',\pi_0'',\pi'')$
will be $1$ if either $p \ne q$ or $r \ne s$.  Otherwise, the distance
is
\be{Mdis2}
1\wedge\left(
  \frac{1}{2}\inf_{\Re^{\mathrm{root,cuts}}}{\mathrm{dis}}
  (\Re^{\mathrm{root,cuts}})
   +\max_i|\sigma_i'-\sigma_i''|
   +\max_j|\tau_j' - \tau_j''| \right),
\ee
where the infimum is over all correspondences 
between $T'$ and $T''$ that contain 
the pairs
$(\rho_{T'},\rho_{T''})$ and
$(x_i',x_i'')$ for $1\le i\le p$.

Equip $\Omega^{\mathrm{fin}}$ with the Borel
$\sigma$-field corresponding to this metric.
For $t \ge 0$, let ${\mathcal F}_t^o$ be the
$\sigma$-field on $\Omega$ generated by the family of maps
from $\Omega$ into  $\Omega^{\mathrm{fin}}$ given by
$(T,\pi_0,\pi) 
\mapsto 
(R_\eta(T), 
\pi_0 \cap (]0,t] \times (R_\eta(T))^o), 
\pi \cap \{(s,x) : x \le s \le t\})$ 
for $\eta>0$.
As usual, set ${\mathcal F}_t^+ := \bigcap_{u>t} {\mathcal F}_u^o$
for $t \ge 0$.
Put ${\mathcal F}^o := \bigvee_{t \ge 0} {\mathcal F}_t^o$.

It is straightforward to establish the following result
from Lemma~\ref{L1} and the construction of $X$
in Subsection \ref{SS:beginconstruct},
and we omit the proof.

\begin{lemma}
\label{L:Xprogmeas}
The map $(t,(T,\pi_0,\pi)) \mapsto X_t(T,\pi_0,\pi)$
from $\R^+ \times \Omega$ into ${\bf T}^{\mathrm{root}}$
is progressively measurable with respect to
the filtration $({\mathcal F}_t^o)_{t \ge 0}$.  (Here, of course,
we are equipping ${\bf T}^{\mathrm{root}}$ with the
Borel $\sigma$-field associated with the metric 
$d_{{\mathrm{GH}}^{\mathrm{root}}}$.)
\end{lemma}

Given $T \in {\bf T}^{\mathrm{root}}$, let ${\bf P}^T$
be the probability measure on $\Omega$ defined by the
following requirements.
\begin{itemize}
\item
The measure ${\bf P}^T$ assigns all of its mass to the
set $\{(T',\pi_0',\pi') \in \Omega : T' = T\}$.
\item
Under ${\bf P}^T$, the random variable $(T',\pi_0',\pi') \mapsto \pi_0'$
is a Poisson point process on the set $\R^{++} \times T^o$ 
with intensity $\lambda \otimes \mu$, where $\mu$ is the length
measure on $T$.
\item
Under ${\bf P}^T$, the random variable $(T',\pi_0',\pi') \mapsto \pi'$
is a Poisson point process on the set
$\{(t,x)\in\R^{++}\times\R^{++} : x \le t\}$ with
intensity $\lambda \otimes \lambda$ restricted to this set.
\item
The random variables $(T',\pi_0',\pi') \mapsto \pi_0'$
and $(T',\pi_0',\pi') \mapsto \pi'$ are independent under ${\bf P}^T$.
\end{itemize}
Of course, the random variable 
$(T',\pi_0',\pi') \mapsto \pi_0'$
takes values in a space of
equivalence classes of countable sets rather than a space of sets
{\em per se}, so, more formally, this random variable
has the law of the image of a Poisson process on an
arbitrary class
representative under the appropriate quotient map.

For $t \ge 0$, $g$ a bounded Borel function on ${\bf T}^{\mathrm{root}}$,
and $T \in {\bf T}^{\mathrm{root}}$, set
\be{defsemigpX}
P_t g(T) := {\bf P}^T[g(X_t)].
\ee
With a slight abuse of notation, let $\tilde R_\eta$ for $\eta > 0$
also denote the map from $\Omega$ into $\Omega$ that sends
$(T,\pi_0,\pi)$ to $(R_\eta(T), \pi_0 \cap 
(\R^{++} \times (R_\eta(T))^o), \pi)$.

Our main construction result is the following.

\begin{theorem}
\label{main}
\begin{itemize}
\item[(i)]
If $T \in {\bf T}^{\mathrm{root}}$ is finite,
then $(X_t)_{t \ge 0}$ under ${\bf P}^T$ is a Markov process that evolves via
the root growth with re-grafting dynamics on finite trees.
\item[(ii)] 
For all $\eta>0$ and $T \in {\bf T}^{\mathrm{root}}$, the law of $(X_t \circ \tilde R_\eta)_{t \ge 0}$
under ${\bf P}^T$ coincides with the law of $(X_t)_{t \ge 0}$ under ${\bf P}^{R_\eta(T)}$.
\item[(iii)]
For all  $T \in {\bf T}^{\mathrm{root}}$, the law of  $(X_t)_{t \ge 0}$ under ${\bf P}^{R_\eta(T)}$
converges as $\eta \downarrow 0$ to that of  $(X_t)_{t \ge 0}$ under ${\bf P}^T$ (in the sense of convergence of laws
on the space of c\`adl\`ag ${\bf T}^{\mathrm{root}}$-valued paths equipped with the
Skorohod topology).
\item[(iv)] For $g \in {\mathrm b}{\mathcal B}({\bf T}^{\mathrm{root}})$, the map
$(t,T) \mapsto P_t g(T)$ is ${\mathcal B}(\R^+) \times {\mathcal B}({\bf T}^{\mathrm{root}})$-measurable.
\item[(v)]  The process $(X_t, {\bf P}^T)$ is strong Markov with respect to the filtration
$({\mathcal F}_t^+)_{t \ge 0}$ and has transition semigroup $(P_t)_{t \ge 0}$.
\end{itemize}
\end{theorem}

\begin{proof}
(i) This is clear from the definition of the root growth and
re-grafting dynamics.

\bigskip\noi
(ii) It is enough to check that the push-forward
of the probability measure ${\bf P}^T$ under the map
$R_\eta: \Omega \rightarrow \Omega$ is the measure 
${\bf P}^{R_\eta(T)}$.   This, however, follows from the observation
that the restriction of length measure on a tree to a subtree is
just length measure on the subtree.

\bigskip\noi
(iii) This is immediate from part (ii), the limiting construction in Subsection
\ref{SS:beginconstruct}, and part (iv) of  Lemma~\ref{L1}.  Indeed, we have
that
\be{uniformcontrol}
\sup_{t \ge 0} d_{\mathrm{GH}^{\mathrm{root}}}(X_t, X_t \circ \tilde R_\eta) 
\le d_{\mathrm{H}}(T, R_\eta(T)) \le \eta.
\ee{}

\bigskip\noi
(iv) By a monotone class argument, it is enough to consider the
case where the test function $g$ is continuous.  It follows from
part (iii) that $P_t g(R_\eta(T))$ converges pointwise to
$P_t g(T)$ as $\eta \downarrow 0$, and it is not difficult to show
using Lemma~\ref{L1}
and part (i) that $(t,T) \mapsto P_t g(R_\eta(T))$ is  
${\mathcal B}(\R^+) \times {\mathcal B}({\bf T}^{\mathrm{root}})$-measurable.
We omit the details, because we will establish an
even stronger result in Proposition \ref{P:Feller}.

\bigskip\noi
(v) By construction and part (ii) of Lemma~\ref{L00}, we have for $t \ge 0$
and $(T,\pi_0,\pi) \in \Omega$ that, as a set, 
$X_t^o(T,\pi_0,\pi)$ is the disjoint union $T^o\amalg ]0,t]$.

Put
\be{deftheta}
\begin{split}
& \theta_t(T,\pi_0,\pi) \\
& \quad := 
\Bigl(X_t(T,\pi_0,\pi), 
\{(s,x) \in \R^{++} \times T^o : (t+s,x ) \in \pi_0\}, \\
& \qquad \{(s,x) \in \R^{++} \times \R^{++} : (t+s, t+x) \in \pi\}\Bigr) \\
& \quad =
\Bigl(X_t(T,\pi_0,\pi), 
\{(s,x) \in \R^{++} \times  X_t^o(T,\pi_0,\pi) : (t+s,x ) \in \pi_0\}, \\
& \qquad \{(s,x) \in \R^{++} \times \R^{++} : (t+s, t+x) \in \pi\}\Bigr). \\
\end{split}
\ee
Thus $\theta_t$ maps $\Omega$ into $\Omega$. 
Note that $X_s \circ \theta_t = X_{s+t}$ and that
$\theta_s \circ \theta_t = \theta_{s+t}$, that is,
the family $(\theta_t)_{t \ge 0}$ is a semigroup.  It is not hard to show
that $(t,(T,\pi_0,\pi)) \mapsto \theta_t(T,\pi_0,\pi)$ 
is jointly measurable, and we leave this to the reader.

Fix $t \ge 0$ and $(T,\pi_0,\pi) \in \Omega$.
Write $\mu'$ for the measure on $T^o\amalg ]0,t]$ that
restricts to length measure on $T^o$ and 
to Lebesgue measure on $]0,t]$.
Write $\mu''$ for the length measure on $X_t^o(T,\pi_0,\pi)$.
The strong Markov property will follow from a standard
strong Markov property for Poisson processes if we can show that
$\mu'=\mu''$.  This equality is clear from the construction if
$T$ is finite:  the tree $X_t(T,\pi_0,\pi)$ is produced from
the tree $T$ and the set $]0,t]$ by a finite number of
dissections and rearrangements.  The equality for general
$T$ follows from the construction and part (iii) of Lemma~\ref{L00}.
\end{proof}

\section{Connection with Aldous's construction of the CRT}
\label{S:Aldousconnect}

Let $({\mathcal{R}}_t)_{t \ge 0}$, $({\mathcal{T}}_t)_{t \ge 0}$,
and $(\tau_n)_{n \in \N}$ be as in the Introduction.
Thus $({\mathcal{T}}_t)_{t \ge 0}$ has the same law
as $(X_t)_{t \ge 0}$ under ${\mathbf{P}}^{T_0}$,
where $T_0$ is the trivial tree.  

\begin{proposition}
\label{P:Aldousconnect}
The two random finite rooted trees
${\mathcal R}_{\tau_n - }$ and ${\mathcal T}_{\tau_n - }$
have  the same distribution for all $n \in \N$.
\end{proposition}

\begin{proof}
Let $R_n$ denote the object obtained by taking the rooted finite tree with
edge-lengths
${\mathcal R}_{\tau_n - }$ and labeling the leaves with $1, \ldots, n$,
in the order that they are added in Aldous's construction.
Let $T_n$ be derived similarly from the rooted finite tree with 
edge-lengths ${\mathcal T}_{\tau_n - }$,
by labeling the leaves with $1, \ldots, n$ in the order that they
appear in the root growth with re-grafting construction.
It will suffice to show that $R_n$ and $T_n$ have the same distribution.
Note that both $R_n$ and $T_n$ are rooted, bifurcating trees with $n$
labeled leaves and edge-lengths.
Such a tree $S_n$ is uniquely specified by its {\em shape},
denoted $\shape(S_n)$,
which is a rooted, bifurcating, leaf-labeled combinatorial tree, 
and by the list of its $(2n-1)$ edge-lengths in a canonical order determined 
by its shape, say
$$
\lengths(S_n) := ( \length(S_n,1), \ldots, \length(S_n, 2n-1) ),
$$
where the edge-lengths are listed in order of traversal of edges by first 
working along the path from the root to leaf $1$, then along the path 
joining that path to leaf $2$, and so on.

Recall that $\tau_n$ is the $n$th point of a Poisson process on $\R^{++}$ 
with rate $t\,dt$.  We construct $R_n$ and $T_n$ on the same
probability space using cuts at points $U_i \tau_i$, $1 \le i \le n-1$,
where $U_1, U_2, \ldots$ is a sequence of independent random variables
uniformly distributed on the interval $]0,1]$ and independent of the 
sequence $(\tau_n)_{n \in \N}$.
Then, by construction, the common collection of edge-lengths of 
$R_n$ and of $T_n$ is the collection of lengths of the $2n-1$ subintervals of
$]0,\tau_n]$
obtained by cutting this interval at the $2n-2$ points
$$
\{ X_i^{(n)}, 1 \le i \le 2n-2 \} := \bigcup_{i = 1}^{n-1} \{ U_i \tau_i , \tau_i \}
$$
where the $X_i^{(n)}$ are indexed to increase in $i$ for each fixed $n$. 
Let $X_0^{(n)}:= 0$ and $X_{2n-1}^{(n)}:= \tau_n$.
Then
\be{lseq}
\length(R_n, i ) = X_i^{(n)} - X_{i-1}^{(n)}, \quad 1 \le i \le 2n  - 1, 
\ee
\be{lseq1}
\length(T_n, i ) = \length(R_n, \sigma_{n,i} ), \quad 1 \le i \le 2n  - 1, 
\ee
for some almost surely unique random indices
$\sigma_{n,i} \in \{1 , \ldots 2n  - 1 \}$ such that
$i \mapsto \sigma_{n,i}$ is almost surely a permutation of 
$\{1 , \ldots 2n  - 1 \}$.
According to \cite[Lemma 21]{Ald93},  the distribution of $R_n$ may
be characterized as follows:
\begin{itemize}
\item[(i)]
the sequence $\lengths(R_n)$ is exchangeable, with the same distribution as the
sequence of lengths of subintervals obtained by cutting $]0,\tau_n]$ at $2n-2$
uniformly chosen points $\{U_i \tau_n: \, 1 \le i \le 2n-2\}$;
\item[(ii)]
$\shape(R_n)$ is uniformly distributed on the set of all
$1 \times 3 \times 5 \times \cdots \times (2 n - 3)$ possible shapes;
\item[(iii)]
$\lengths(R_n)$ and $\shape(R_n)$ are independent.
\end{itemize}
In view of this characterization and \re{lseq1},
to show that $T_n$ has the same distribution as $R_n$ it is
enough to show that 
\begin{itemize}
\item[(a)]
the random permutation $\{i \mapsto \sigma_{n,i}: \, 1 \le i \le 2n - 1 \}$ is a
function
of $\shape(T_n)$;
\item[(b)]
$\shape ( T_n) = \Psi_n( \shape(R_n) )$
for some bijective map $\Psi_n$ from the set of all possible shapes to itself.
\end{itemize}
This is trivial for $n = 1$, so we assume below that $n \ge 2$.
Before proving (a) and (b), we recall that (ii) above involves a 
natural bijection
\be{bij}
(I_1, \ldots, I_{n-1} ) \leftrightarrow 
\shape(R_n) 
\ee
where $I_{n-1} \in \{1, \ldots, 2n - 3\}$ is the unique $i$ such that
$U_{n-1}\tau_{n-1} \in (X_{i-1}^{(n-1)}, X_i^{(n-1)})$. Hence $I_{n-1}$ is 
the index in the canonical ordering of edges of $R_{n-1}$
of the edge that is cut in the transformation from $R_{n-1}$ to $R_n$ 
by attachment of an additional edge,
of length $\tau_n - \tau_{n-1}$, connecting the cut-point to leaf $n$.
Thus (ii) and (iii) above correspond via \re{bij} to the facts 
that $I_1, \ldots, I_{n-1}$ are independent and uniformly distributed
over their ranges, and independent of $\lengths(R_n)$. These facts can be 
checked directly from the construction of $(R_n)_{n \in \N}$
from $(\tau_n)_{n \in \N}$ and $(U_n)_{n \in \N}$ 
using standard facts about uniform order statistics.

Now (a) and (b) follow from  \re{bij} and
another bijection
\be{bij1}
(I_1, \ldots, I_{n-1} ) 
\leftrightarrow 
\shape(T_n)  
\ee
where each possible value $i$ of $I_{m}$ is identified with edge
$\sigma_{m,i}$ in the canonical ordering of edges of $T_{m}$. This is
the edge of $T_m$ whose length equals $\length(R_m,i)$.
The bijection \re{bij1},
and the fact that $\sigma_{n,i}$ depends only on $\shape(T_n)$,
will now be established by induction on $n \ge 2$.
For $n = 2$ the claim is obvious.
Suppose for some $n \ge 3$ that the correspondence between $(I_1, \ldots, I_{n-2} )$
and $\shape(T_{n-1})$ has been established, and that 
the length of edge $\sigma_{n-1,i}$ in the canonical ordering of edges of
$T_{n-1}$ is equals the length of the $i$th edge in the canonical ordering
of edges of $R_{n-1}$, for some $\sigma_{n-1,i}$ which is a function of
$i$ and $\shape(T_{n-1})$.
According to the construction of $T_n$, 
if $I_{n-1} = i$ then $T_n$ is derived from $T_{n-1}$ by splitting $T_{n-1}$ 
into two branches at some point along edge $\sigma_{n-1,i}$ in the
canonical ordering of the edges of 
$T_{n-1}$, and forming a new tree from the two 
branches and an extra segment of length $\tau_n - \tau_{n-1}$.  
Clearly, $\shape(T_n)$ is determined by $\shape(T_{n-1})$ and
$I_{n-1}$, and in the canonical ordering of the edge-lengths of $T_n$ the
length of the $i$th edge equals the length of the edge $\sigma_{n,i}$ of 
$R_{n}$, for some $\sigma_{n,i}$ which is a function of $\shape(T_{n-1})$ and
$I_{n-1}$, and hence a function of $\shape(T_n)$.  To complete the proof,
it is enough
by the inductive hypothesis to show that the map
$$
( \shape(T_{n-1}), I_{n-1} ) \rightarrow \shape(T_n)  
$$
just described is invertible. But $\shape(T_{n-1})$ and $I_{n-1}$ can be recovered
from $\shape(T_n)$ by the following sequence of moves:
\begin{itemize}
\item delete the edge attached to the root of $\shape(T_n)$ 
\item split the remaining tree into its two branches leading away from the 
internal node to which the deleted edge was attached;
\item re-attach the bottom end of the branch not containing leaf $n$ to
leaf $n$ on the other branch, joining the two incident edges to form a
single edge; 
\item the resulting shape is $\shape(T_{n-1})$, and $I_{n-1}$ is the
index such that the joined edge in $\shape(T_{n-1})$ is the edge
$\sigma_{n-1,I_{n-1}}$ in the canonical ordering of edges on
$\shape(T_{n-1})$.
\end{itemize}
\end{proof}

\section{Recurrence and convergence to stationarity}
\label{S:recurstat}

\begin{lemma}
\label{L:coupletotrivial}
For any $(T,d,\rho) \in {\mathbf T}^{\mathrm{root}}$
we can build on the same probability space two 
${\mathbf T}^{\mathrm{root}}$-valued processes $X'$ and $X''$ such that:
\begin{itemize}
\item
$X'$ has the law of $X$ under ${\mathbf P}^{T_0}$, where $T_0$
is the trivial tree  consisting of just the root, 
\item
$X''$ has the law of $X$ under ${\mathbf P}^T$, 
\item
for all $t \ge 0$,
\be{keepclose}
d_{\mathrm{GH}^{\mathrm{root}}}(X_t',X_t'') 
\le d_{\mathrm{GH}^{\mathrm{root}}}(T_0,T)
= \sup\{d(\rho,x) : x \in T\}
\ee
\item
\be{getcloser}
\lim_{t \rightarrow \infty} d_{\mathrm{GH}^{\mathrm{root}}}(X_t',X_t'') = 0, \quad \text{almost surely}.
\ee
\end{itemize}
\end{lemma}

\begin{proof}
The proof  follows almost immediately from
 construction of $X$ in Section \ref{S:construct}
and Lemma~\ref{L0}.  The only point requiring some comment is (\ref{getcloser}).
For that it will be enough to show for any $\varepsilon > 0$
that for ${\mathbf P}^T$-a.e.
$(T,\pi_0,\pi) \in \Omega$ there exists $t>0$ such that
the projection of
$\pi_0\cap (]0,t] \times T^o)$ onto $T$ is an $\varepsilon$-net for
$T$.

Note that the projection of $\pi_0\cap (]0,t] \times T^o)$ onto $T$
is a Poisson process under  ${\mathbf P}^T$ with intensity
$t \mu$, where $\mu$ is the length measure on $T$.
Moreover, $T$ can be covered by a finite collection of
$\varepsilon$-balls, each with positive $\mu$-measure. 
Therefore, the ${\mathbf P}^T$-probability of the set of
$(T,\pi_0,\pi) \in \Omega$ such that 
the projection of
$\pi_0\cap (]0,t] \times T^o)$ onto $T$ is an $\varepsilon$-net for
$T$ increases as $t \rightarrow \infty$ to $1$.
\end{proof}

\begin{proposition}
\label{P:convtoCRT}
For any $T \in {\mathbf T}^{\mathrm{root}}$, the law of
$X_t$ under ${\mathbf P}^T$ converges weakly to that of the
Brownian CRT as $t \rightarrow \infty$.
\end{proposition}

\begin{proof}
It suffices by
 Lemma~\ref{L:coupletotrivial} to consider the case where
$T$ is the trivial tree.

We saw in the Proposition~\ref{P:Aldousconnect} that,
in the notation of the Introduction, ${\mathcal T}_{\tau_n-}$ has the same
distribution as ${\mathcal R}_{\tau_n-}$.  Moreover, we recalled
in the Introduction that
${\mathcal R}_t$ converges in distribution to the continuum
random tree as $t \rightarrow \infty$ if we use Aldous's
metric on trees that comes from thinking of them as closed subsets
of $\ell^1$ with the root at the origin and equipped with the
Hausdorff distance.  By construction, $({\mathcal T}_t)_{t \ge 0}$
has the root growth with re-grafting dynamics started at the trivial
tree.  Clearly, the rooted Gromov--Hausdorff distance between
${\mathcal T}_t$ and ${\mathcal T}_{\tau_{n+1}-}$ is at most
$\tau_{n+1} - \tau_n$ for $\tau_n \le t < \tau_{n+1}$.  It remains
to observe that $\tau_{n+1} - \tau_n \rightarrow 0$ in probability
as $n \rightarrow \infty$.
\end{proof}

\begin{proposition} 
\label{P:toprecur}
Consider a non-empty open set $U \subseteq {\mathbf T}^{\mathrm{root}}$.
For each $T \in {\mathbf T}^{\mathrm{root}}$,
\be{recurstatement}
{\mathbf P}^T\{
\text{for all $s \ge 0$, there exists $t > s$ such that  $X_t \in U$}\} = 1.
\ee
\end{proposition}

\begin{proof}
It is straightforward, but notationally
rather tedious, to show that if 
$B' \subseteq {\mathbf T}^{\mathrm{root}}$ is any ball
and $T_0$ is the trivial tree,
then 
\be{supportproperty}
{\mathbf P}^{T_0}\{X_t \in B'\} > 0
\ee
for all $t$ sufficiently large.
Thus, for any ball 
$B' \subseteq {\mathbf T}^{\mathrm{root}}$ there is,
by Lemma~\ref{L:coupletotrivial},
a ball $B'' \subseteq {\mathbf T}^{\mathrm{root}}$
containing the trivial tree
such that
\be{localsupportprop}
\inf_{T \in B''} {\mathbf P}^T\{X_t \in B'\} > 0
\ee
for each $t$ sufficiently large.

By a standard application of the Markov property,
it therefore suffices to show for each $T \in {\mathbf T}^{\mathrm{root}}$
and each ball $B''$ around the trivial tree that
\be{musthitball}
{\mathbf P}^T\{\text{there exists $t>0$ such that $X_t \in B''$}\}=1.
\ee

By another standard application of the Markov property, equation
(\ref{musthitball}) will follow if we can show that there
is a constant $p>0$ depending
on $B''$ such that for any $T \in {\mathbf T}^{\mathrm{root}}$
\begin{equation*}
\liminf_{t \rightarrow \infty} {\mathbf P}^T\{X_t \in B''\} > p.
\end{equation*}
This, however, follows from Proposition \ref{P:convtoCRT} and the observation
that for any $\varepsilon > 0$ 
the law of the Brownian CRT assigns positive mass to the set of
trees with height less than $\varepsilon$ (which is just the observation
that the law of the Brownian excursion assigns positive mass to the
set of excursion paths with maximum less that $\varepsilon/2$).
\end{proof}

\begin{proposition}
\label{P:CRTstationary}
The law of the Brownian CRT is the unique stationary
distribution for $X$.  That is, if $\xi$ is
the law of the CRT, then 
$\int \xi(dT) P_t f(T) = \int \xi(dT) f(T)$
for all $t \ge 0$ and 
$f \in {\mathrm b}{\mathcal B}({\mathbf T}^{\mathrm{root}})$,
and $\xi$ is the unique probability measure on ${\mathbf T}^{\mathrm{root}}$
with this property.
\end{proposition}

\begin{proof}
This is a standard argument given Proposition~\ref{P:convtoCRT}
and the Feller property for the semigroup $(P_t)_{t \ge 0}$
established in Proposition~\ref{P:Feller}, but we include
the details for completeness.

Consider a test function 
$f: {\mathbf T}^{\mathrm{root}} \rightarrow \R$ 
that is continuous and bounded.
By Proposition~\ref{P:Feller} below, the
function $P_t f$ is also continuous and bounded for each
$t \ge 0$.
Therefore, by Proposition~\ref{P:convtoCRT}, 
\be{}
\begin{split}
\int \xi(dT) f(T) & = \lim_{s \rightarrow \infty} \int \xi(dT) P_s f(T)
= \lim_{s \rightarrow \infty} \int \xi(dT) P_{s+t} f(T) \\
& = \lim_{s \rightarrow \infty} \int \xi(dT) P_s (P_t f)(T)
= \int \xi(dT) P_t f(T) \\
\end{split}
\ee
for each $t \ge 0$, and hence $\xi$ is stationary.  Moreover,
if $\zeta$ is a stationary measure, then
\be{}
\begin{split}
\int \zeta(dT) f(T) & = \int \zeta(dT) P_t f(T) \\
& \rightarrow \int \zeta(dT) \left(\int \xi(dT) f(T) \right)
= \int \xi(dT) f(T), \\
\end{split}
\ee
and $\zeta=\xi$, as claimed. 
\end{proof}

\section{Feller property}
\label{S:Feller}

The following result says that the law of $X_t$ under
${\mathbf P}^T$ is weakly continuous in 
the initial value $T$ for each
$t \ge 0$.  This property is sometimes referred to
as the Feller property of the semigroup $(P_t)_{t \ge 0}$,
although this terminology is often restricted to the case
of a locally compact state space and transition operators that 
map the space of continuous functions that vanish at infinity
into itself.  A standard consequence of this result
is that the law of the process $(X_t)_{t \ge 0}$ is weakly
continuous in the initial value (when the space of
c\`adl\`ag ${\mathbf T}^{\mathrm{root}}$-valued paths
is equipped with the Skorohod topology).

\begin{proposition}
\label{P:Feller}
If the function $f: {\mathbf T}^{\mathrm{root}} \rightarrow \R$ 
is continuous and bounded, then the
function $P_t f$ is also continuous and bounded for each
$t \ge 0$.
\end{proposition}

We will prove the proposition by a coupling argument that,
{\em inter alia}, builds
processes with the law of $X$ under ${\mathbf P}^T$
for two different finite values of $T$ on the same
probability space.  The key to constructing such a coupling
is the following pair of lemmas.

We require the following notion. A {\em rooted combinatorial tree}
is just a connected, acyclic graph
with one vertex designated as the root.  Equivalently,
we can think of
a rooted combinatorial tree as a finite rooted tree
in which all edges have length one.  Thus any finite rooted
tree is associated with a unique rooted combinatorial tree by
changing all the edge lengths to one, and any two finite
rooted trees with the same topology are associated with the
same rooted combinatorial tree.  If $U$ and $V$ are
two rooted combinatorial trees with leaves labeled
by $(x_1, \ldots, x_n)$ and $(y_1, \ldots, y_n)$, then
we say that $U$ and $V$ are isomorphic if there
exists a graph isomorphism between $U$ and $V$
that maps the root of $U$ to the root of $V$
and $x_i$ to $y_i$ for $1 \le i \le n$.

\begin{lemma}\label{L:Tdoubprimeexist}
Let $(T,\rho)$ be a finite rooted trees with leaves $\{x_1,\ldots,x_n\}=T\setminus T^o$.
(recall the definition of the skeleton $T^o$ from (\ref{sce})).  Write
$\eta$ for the minimum of the (strictly positive) edge lengths in $T$. Suppose that
$(T',\rho')$ is another finite rooted tree with 
$d_{{\mathrm{GH}}^{\mathrm{root}}}((T^\prime,\rho^\prime),(T,\rho))<\delta<\frac{\eta}{16}$.
Then there exists a subtree $(T^{\prime\prime},\rho^\prime)\preceq^{\mathrm{root}}(T^{\prime},\rho^\prime)$
and a map $\bar f:T\to T^{\prime\prime}$ such that:
\begin{itemize}
\item[(i)] $\bar f(\rho) = \rho'$,
\item[(ii)] $T^{\prime\prime}$ is spanned by $\{\bar f(x_1),\ldots, \bar f(x_n),\rho^{\prime}\}$,
\item[(iii)] $d_{\mathrm{H}}(T',T'') < 3 \delta$,
\item[(iv)] ${\mathrm{dis}}(\bar f) < 8 \delta$,
\item[(v)] $T''$ has leaves $\{\bar f(x_1),\ldots, \bar f(x_n)\}$,
\item[(vi)] 
by possibly deleting some internal edges
from the rooted combinatorial tree associated to $T''$
with leaves labeled by
$(\bar f(x_1), \ldots, \bar f(x_n))$,
one can obtain a 
leaf-labeled rooted combinatorial tree that is isomorphic to the rooted
combinatorial rooted tree 
associated to $T$ with leaves labeled by
$(x_1, \ldots, x_n)$. 
\end{itemize}
\end{lemma}

\begin{proof}
We have from (\ref{GHroot})
that there is a correspondence ${\Re}^{\mathrm{root}}$ containing 
$(\rho,\rho^\prime)$ between $T$ and $T'$ such that ${\mathrm{dis}}({\Re}^{\mathrm{root}})<2\delta$.
For $x\in T\setminus\{\rho\}$, choose $f(x)\in T^\prime$ such that
$(x,f(x))\in{\Re}^{\mathrm{root}}$, and put $f(\rho):=\rho^\prime$. 
Set $T^{\prime\prime}$ to be the subtree of
$T'$  spanned by $\{f(x_1),\ldots, f(x_n),\rho^{\prime}\}$.
For $x \in T$ define $\bar f(x) \in T''$ to be the point in $T''$
that has minimum distance to $f(x)$.  In particular,
$\bar f(\rho) = f(\rho) = \rho'$ and $\bar f(x_i) = f(x_i)$ for all $i$,
so that (i) and (ii) hold.

For $x^\prime\in T^\prime\setminus\{\rho^\prime\}$
choose $g(x^\prime)$ such that
$(g(x^\prime),x^\prime)\in{\Re}^{\mathrm{root}}$
and put $g(\rho^\prime):=\rho$.
Then $f(T)$ and $g(T^\prime)$ are $2\delta$-nets for $T^\prime$ and $T$, 
respectively, and ${\mathrm{dis}}(f)\vee{\mathrm{dis}}(g)<2\delta$. 
For each $i\in\{1,\ldots,n\}$ and $y^\prime\in T^\prime$ we have 
$(\rho,\rho^\prime),(x_i,f(x_i)),(g(y^\prime),y^\prime)\in\Re^{\mathrm{root}}$.
Hence $d_{T^\prime}(\rho^\prime,y^\prime)<d_{T}(\rho,g(y^\prime))+2\delta$,
and $d_{T^\prime}(y^\prime,f(x_i))<d_{T}(g(y^\prime),x_i)+2\delta$.
Now fix $y^\prime\in T^\prime$, and choose 
$i\in\{1,\ldots,n\}$ such that $g(y^\prime)\in[\rho,x_i]$. Then 
\be{est16}
\begin{aligned}
   d_{T^\prime}(\rho^\prime,f(x_i))+
   2d_{\mathrm H}\left(\{y^\prime\},[\rho^\prime,f(x_i)]\right)
 &=
   d_{T^\prime}(\rho^\prime,y^\prime)+
   d_{T^\prime}(y^\prime,f(x_i))
  \\
 &<
   d_T(\rho,x_i)+4\delta
   \\
 & <
   d_{T'}(\rho',f(x_i)) + 2 \delta + 4 \delta,
\end{aligned}
\ee
and hence $d_{\mathrm H}(\{y^\prime\},T^{\prime\prime})< 3\delta$. Thus (iii) holds.

For $x,y\in T$, 
\be{est23}
\begin{aligned}
   |d_{T}&(x,y)-d_{T^{\prime\prime}}(\bar{f}(x),\bar{f}(y))|
  \\[1mm]
 &\le
   |d_{T}(x,y)-d_{T^{\prime}}(f(x),f(y))|
   +d_{T^{\prime}}(\bar{f}(x),f(x))+d_{T^{\prime}}(\bar{f}(y),f(y))
  \\
 &\le
   {\mathrm{dis}}(f)+2d_{H}(T^{\prime},T^{\prime\prime})<8\delta,
\end{aligned}
\ee
and (iv) holds.

In order to establish (v), it suffices to observe for $1 \le i \ne j \le n$
that, by part (iv),
\be{stricttriangle}
\begin{split}
& d_{T''}(\bar f(x_i), \bar f(x_j)) + d_{T''}(\bar f(x_j), \rho') - d_{T''}(\bar f(x_i), \rho') \\
& \quad \ge
d_{T}(x_i, x_j) + d_{T}(x_j, \rho) - d_{T}(x_i, \rho) - 3 {\mathrm{dis}}(\bar f) \\
& \quad >
2 \eta - 24 \delta \\
& \quad > 0. \\
\end{split}
\ee

Similarly, part (vi) follows from part (iv) and the observations in Subsection \ref{fourp} about
re-constructing tree shapes from distances between the points in subsets of size
four drawn from the leaves and the root of $T''$ once we observe the inequality
$\frac{1}{2} 4 {\mathrm{dis}}(\bar f) < 16 \delta < \eta$.
\end{proof}

\begin{lemma}\label{L4} 
Let $(T,\rho)$ be a finite rooted tree and
  $\varepsilon>0$. There exists $\delta>0$ depending on $T$ and $\varepsilon$ such
  that if $(T^\prime,\rho^\prime)$ is
  a finite rooted tree with 
  $d_{{\mathrm{GH}}^{\mathrm{root}}}((T^\prime,\rho^\prime),(T,\rho))<\delta$, then there
  exist subtrees $(S,\rho)\preceq^{\mathrm{root}}T$ and
  $(S^\prime,\rho^\prime)\preceq^{\mathrm{root}}T^\prime$ 
  for which:
\begin{itemize}
\item[(i)] $d_{\mathrm H}(S,T)<\varepsilon$ and
  $d_{\mathrm H}(S^\prime,T^\prime)<\varepsilon$, 
\item[(ii)] $S$ and $S^\prime$ have the same total length,
\item[(iii)] there is a bijective measurable map $\psi:S\to S^\prime$
that preserves length measure and has distortion at most $\varepsilon$,
\item[(iv)] the length measure of the set of points $a \in S$ such that
$\{b' \in S': \psi(a) \le b'\} \ne \psi(\{b \in S: a \le b\})$
(that is, the set of points $a$
such that the subtree above $\psi(a)$ is not the image 
under $\psi$ of the subtree above $a$) is less than $\varepsilon$.
\end{itemize}
\end{lemma}

\begin{proof}
As in  Lemma~\ref{L:Tdoubprimeexist}, denote by $\eta$ the minimum of
the (strictly positive) edge lengths of $T$.
Let $(T',\rho')$ be a finite rooted tree with 
\begin{equation}
d_{{\mathrm{GH}}^{\mathrm{root}}}((T^\prime,\rho^\prime),(T,\rho))<\delta< \frac{\eta}{16},
\end{equation}
where $\delta$ depending on $T$ and $\varepsilon$ will be chosen later.
Set $(T'',\rho')$ and $\bar f$ to be a subtree of $T'$ and a function
from $T$ to $T''$ whose existence is guaranteed by  Lemma~\ref{L:Tdoubprimeexist} for this
choice of $\delta$.  Let $\{x_1, \ldots x_n\}$
denote the leaves of $T$ and write  $x_i^\prime:=f(x_i)=\bar{f}(x_i)$ for $i=1,\ldots,n$.

Define inductively subtrees $S_1,\ldots,S_n$ of $T$ (all with root $\rho$) and
$S_1^\prime,\ldots,S_n^\prime$ of $T^{\prime\prime}\subseteq T^\prime$ 
(all with root $\rho^\prime$) as follows. 
Set $S_1:=[\rho,y_1]$
and $S_1^\prime:=[\rho,y_1^\prime]$, where $y_1$ and $y_1^\prime$ are
the unique points on the arcs $[\rho,x_1]$ and
$[\rho^\prime,x_1^\prime]$,
respectively, such that
\be{est17}
   d_T(\rho,y_1)=d_{T^\prime}(\rho^\prime,y_1^\prime)=
   d_T(\rho,x_1)\wedge d_{T^\prime}(\rho^\prime,x_1^\prime).
\ee
Suppose that $S_1,\ldots,S_m$ and $S_1^\prime,\ldots,S_m^\prime$ have been
defined. Let $z_{m+1}$ and $z_{m+1}^\prime$ be the points on $S_m$ and 
$S_m^\prime$ closest to $x_{m+1}$ and $x_{m+1}^\prime$. Put
$S_{m+1}:=S_m\cup]z_{m+1},y_{m+1}]$ and 
$S_{m+1}^\prime:=S_m^\prime\cup]z_{m+1}^\prime,y_{m+1}^\prime]$, where 
$y_{m+1}$ and $y_{m+1}^\prime$ are the unique points on the arcs 
$]z_{m+1},x_{m+1}]$ and $]z_{m+1}^\prime,x_{m+1}^\prime]$,
respectively, such that 
\be{est18}
\begin{aligned}
   d_{T^\prime}(z_{m+1},y_{m+1})
 &=
   d_{T^\prime}(z_{m+1}^\prime,y_{m+1}^\prime)
  \\
 &=
   d_T(z_{m+1},x_{m+1})\wedge
 d_{T^\prime}(z_{m+1}^\prime,x_{m+1}^\prime).
\end{aligned}
\ee
Set $S:=S_n$ and $S^\prime:=S_n^\prime$.

Put $z_1:=\rho$, and $z_1^\prime:=\rho^\prime$.
By construction, the arcs $]z_k,y_k]$, $1 \le k \le n$, are disjoint
and their union is $S \setminus \{\rho\}$.
Similarly, the arcs $]z_k^\prime,y_k^\prime]$ 
are disjoint and their union  is $S' \setminus \{\rho'\}$.
Moreover, the arcs $]z_k,y_k]$ and $]z_k',y_k']$ have the same
length  (in particular, $S$ and $S'$
have the same length and part (ii) holds).
We may therefore define a measure-preserving bijection $\psi$ between 
$S$ and $S'$ by setting $\psi(\rho) = \rho'$ and letting the
restriction of $\psi$ to each arc  $]z_k,y_k]$ be the obvious
length preserving bijection onto $]z_k',y_k']$.
More precisely, if $a \in ]z_k,y_k]$, then
$\psi(a)$ is the uniquely determined
point on $]z_k^\prime,y_k^\prime]$ such that $d_{S^\prime}(z_k^\prime,\psi(a))=
d_S(z_k,a)$.

We next estimate the distortion of $\psi$ to establish part (iii).
We first claim that for $a,b \in S$,
\be{est00}
   |d_S(a,b)-d_{S^\prime}(\psi(a),\psi(b))|
 \le
   5 \gamma,
\ee
where
\be{}
\gamma := \max_{1 \le k,m \le n} |d_S(y_k,y_m)-d_{S^\prime}(y_k^\prime,y_m^\prime)| \vee
\max_{1 \le k \le n} |d_S(y_k,\rho)-d_{S^\prime}(y_k^\prime,\rho')|.
\ee

To see (\ref{est00}), consider $a,b \in S \setminus \{\rho\}$ with $a \in ]z_k, y_k]$
and $b \in ]z_m, y_m]$ where $k \ne m$.
(The case where $a=\rho$ or $b=\rho$ holds ``by continuity'' and is left
to the reader.)  Without loss of generality, assume that $k<m$, so that
$y_k \wedge y_m \le z_m < b \le y_m$ in the partial order on $S$
and
$y_k' \wedge y_m' \le z_m' < \psi(b) \le y_m'$ in the partial order on $S'$.
Note that $y_k \wedge y_m$ and $z_k$ are comparable in the partial order,
as are $y_k' \wedge y_m'$ and $z_k'$.  Moreover, by part (vi) of  Lemma
\ref{L:Tdoubprimeexist}, $y_k \wedge y_m \le z_k$ 
if and only if  $y_k' \wedge y_m' \le z_k'$.
We then have to consider four cases depending on the relative positions
of $y_k \wedge y_m,  a$ and  $y_k' \wedge y_m', \psi(a)$. 

\bi\noi
{\em Case I:} $y_k \wedge y_m  < a \le y_k$ and $y_k' \wedge y_m'  < \psi(a) \le y_k'$.\\
We have
\be{}
d_S(y_k,y_m) = d_S(y_k,a) + d_S(a,b) + d_S(b,y_m)
\ee
and
\be{}
d_{S'}(y_k',y_m') = d_{S'}(y_k',\psi(a)) + d_{S'}(\psi(a),\psi(b)) + d_{S'}(\psi(b),y_m').
\ee
By construction,
\be{}
d_S(y_k,a) = d_{S'}(y_k',\psi(a))
\ee
and
\be{}
d_S(b,y_m) = d_{S'}(\psi(b),y_m').
\ee
Hence
\be{}
|d_S(a,b) - d_{S'}(\psi(a),\psi(b))| = |d_S(y_k,y_m) - d_{S'}(y_k',y_m')| \le \gamma.
\ee

\bi\noi
{\em Case II:} $y_k \wedge y_m  < a \le y_k$ and $\psi(a) \le y_k' \wedge y_m'  < y_k'$.\\
Note that in this case $z_k \le y_k \wedge y_m$.
We again have
\be{}
d_S(y_k,y_m) = d_S(y_k,a) + d_S(a,b) + d_S(b,y_m),
\ee
but now
\be{}
\begin{split}
d_{S'}(y_k',y_m') & = d_{S'}(y_k',\psi(a)) + d_{S'}(\psi(a),\psi(b)) + d_{S'}(\psi(b),y_m')\\
& \quad - 2 d_{S'}(\psi(a), y_k'\wedge y_m'). \\
\end{split}
\ee
Let $y_\ell$ be such that $z_k = y_\ell \wedge y_k = y_\ell \wedge y_m$ and hence
$z_k' = y_\ell' \wedge y_k' = y_\ell' \wedge y_m'$.
Observe from Subsection \ref{fourp} that 
\be{changerelbdII}
\begin{split}
& d_{S'}(\psi(a), y_k'\wedge y_m') \\
& \quad =
\frac{1}{2}(d_{S'}(y_\ell',y_m') + d_{S'}(y_k',\rho') - d_{S'}(y_k',y_m') - d_{S'}(y_\ell',\rho'))
- d_{S'}(z_k', \psi(a)) \\
& \quad \le
\frac{1}{2}(d_{S}(y_\ell,y_m) + d_{S}(y_k,\rho) - d_{S}(y_k,y_m) - d_{S}(y_\ell,\rho))
+ \frac{1}{2} 4 \gamma
- d_S(z_k, a) \\
& \quad =
d_S(z_k, y_k \wedge y_m) - d_S(z_k,a) + 2 \gamma \\
& \quad \le 2 \gamma, \\
\end{split}
\ee
and hence
\be{}
|d_S(a,b) - d_{S'}(\psi(a),\psi(b))|  \le 5 \gamma.
\ee

\bi\noi
{\em Case III:} $a \le  y_k \wedge y_m  <  y_k$ and $y_k' \wedge y_m' \le \psi(a) < y_k'$.\\
Note that in this case, $z_k' \le y_k' \wedge y_m'$.
This case is similar to Case II, but we record some of the
details for use later in the proof of part (iv).  Letting the index $\ell$ be
as in Case II, we have
\be{}
\begin{split}
d_S(y_k,y_m) & = d_S(y_k,a) + d_S(a,b) + d_S(b,y_m) \\
& \quad - 2 d_S(a, y_k \wedge y_m) \\
\end{split}
\ee
and
\be{}
d_{S'}(y_k',y_m')  = d_{S'}(y_k',\psi(a)) + d_{S'}(\psi(a),\psi(b)) + d_{S'}(\psi(b),y_m').
\ee
We have
\be{changerelbdIII}
\begin{split}
& d_{S}(a, y_k\wedge y_m) \\
& \quad =
\frac{1}{2}(d_{S}(y_\ell,y_m) + d_{S}(y_k,\rho) - d_{S}(y_k,y_m) - d_{S}(y_\ell,\rho))
   - d_S(z_k, a) \\
& \quad \le
\frac{1}{2}(d_{S'}(y_\ell',y_m') + d_{S'}(y_k',\rho') - d_{S'}(y_k',y_m') - d_{S'}(y_\ell',\rho'))
+ \frac{1}{2} 4 \gamma \\
& \qquad - d_{S'}(z_k', \psi(a)) \\
& \quad =
d_{S'}(z_k', y_k' \wedge y_m') - d_{S'}(z_k',\psi(a)) + 2 \gamma \\
& \quad \le 2 \gamma, \\
\end{split}
\ee
and hence
\be{}
|d_S(a,b) - d_{S'}(\psi(a),\psi(b))|  \le 5 \gamma.
\ee

\bi\noi
{\em Case IV:} $a \le  y_k \wedge y_m  <  y_k$ and $\psi(a) \le  y_k' \wedge y_m' < y_k'$.\\
Letting the index $\ell$ be as in Case II, we have
\be{}
d_S(z_k,y_m) = d_S(z_k,a) + d_S(a,b) + d_S(b,y_m)
\ee{}
and
\be{}
d_{S'}(z_k',y_m') = d_{S'}(z_k',\psi(a)) + d_{S'}(\psi(a),\psi(b)) + d_{S'}(\psi(b),y_m').
\ee
Hence, from Subsection \ref{fourp},
\be{}
\begin{split}
& d_S(a,b) - d_{S'}(\psi(a),\psi(b)) \\
& \quad = d_S(z_k,y_m) - d_{S'}(z_k',y_m') \\
& \quad = d_S(y_\ell \wedge y_m, y_m) - d_{S'}(y_\ell' \wedge y_m', y_m') \\
& \quad = \frac{1}{2}(d_S(y_m,\rho) + d_S(y_\ell,y_m) - d_S(y_\ell,\rho)) \\
& \qquad - \frac{1}{2}(d_{S'}(y_m',\rho') + d_{S'}(y_\ell',y_m') - d_{S'}(y_\ell',\rho')). \\
\end{split}
\ee
Thus
\be{}
|d_S(a,b) - d_{S'}(\psi(a),\psi(b))|  \le \frac{3}{2} \gamma.
\ee

Combining Cases I--IV, we see that (\ref{est00}) holds.  
We thus require an estimate of $\gamma$ to complete the estimation
of the distortion of $\psi$ .  
Clearly,
\be{compareys}
\begin{split}
|d_S(y_k,y_m)-d_{S^\prime}(y_k^\prime,y_m^\prime)|
& \le
|d_T(x_k,x_m)-d_{T''}(x_k^\prime,x_m^\prime)| \\
& \quad + d_T(y_k,x_k) + d_T(y_m,x_m) \\
& \quad + d_{T''}(y_k',x_k') + d_{T''}(y_m',x_m'). \\
\end{split}
\ee
By (\ref{est17}),  
\be{}
\begin{split}
d_T(y_1,x_1)\vee d_{T^{\prime\prime}}(y_1',x_1')
& =|d_T(\rho,x_1)- d_{T^{\prime\prime}}(\rho^\prime,x_1^\prime)| \\
& \le{\mathrm{dis}}(\bar{f}) \\
& < 8 \delta. \\
\end{split}
\ee
For $2 \le k \le n$ there exists by construction an index
$i\in\{1,2,\ldots,k-1\}$ such that $z_k\in[z_i,y_i]$ and 
$z_k^\prime\in[z_i^\prime,y_i^\prime]$. Applying the observations
of Subsection \ref{fourp},
\be{est31}
\begin{aligned}
   d_T&(y_k,x_k)\vee d_{T^{\prime\prime}}(y_k^\prime,x_k^\prime)
  \\
 &=
   |d_T(y_k,x_k)-d_{T^{\prime\prime}}(y_k^\prime,x_k^\prime)|
  \\
 &=
   |d_T(z_k,x_k)-d_{T^{\prime\prime}}(z_k^\prime,x_k^\prime)|
  \\
 &\le 
   \frac{1}{2}\{|d_T(x_i,x_k)-d_{T^{\prime\prime}}(x_i^\prime,x_k^\prime)|
  \\
 &\quad +
   |d_T(\rho,x_i)-d_{T^{\prime\prime}}(\rho^\prime,x_i^\prime)|+
   |d_T(\rho,x_k)-d_{T^{\prime\prime}}(\rho^\prime,x_k^\prime)|\}  
  \\
 &\le
   \frac{3}{2}{\mathrm{dis}}(\bar{f})\\
 &\le 12 \delta.
\end{aligned}
\ee
Thus, from (\ref{compareys}),
\be{}
|d_S(y_k,y_m)-d_{S^\prime}(y_k^\prime,y_m^\prime)| < (8 + 4 \times 12)\delta = 56 \delta.
\ee
A similar argument shows that $|d_S(y_k,\rho)-d_{S^\prime}(y_k^\prime,\rho')| < (8 + 2 \times 12)\delta =  32\delta$,
and hence $\gamma < 56 \delta$.
Substituting into (\ref{est00}) gives
\be{}
{\mathrm{dis}}(\psi) \le 5 \gamma < (5 \times 56) \delta = 280 \delta.
\ee

Moving to part (i), apply (\ref{est31}) to obtain
\be{est34}
d_{\mathrm H}(S,T) \le \max_{1 \le i \le n} d(y_k,x_k) \le \gamma < 56 \delta
\ee 
and, by similar arguments, 
\be{est35}
       d_{\mathrm{H}}(S',T') 
   \le d_{\mathrm{H}}(S',T'')
      +d_{\mathrm{H}}(T'',T')
   < 59 \delta.
\ee

Finally, we consider part (iv).  Suppose that
$a \in S$ is such that the subtree of $S'$ above $\psi(a)$ is not the image 
under $\psi$ of the subtree of $S$ above $a$.
Let $k$ be the unique index such that $a \in ]z_k, y_k]$ (and hence
$\psi(a) \in ]z_k', y_k']$).
It follows from the construction of $\psi$ that
there must exist an index $\ell$
such that either
$z_k < a \le z_\ell$ and $z_k' < z_\ell' \le \psi(a)$
or
$z_k < z_\ell \le a$ and $z_k' < \psi(a) \le z_\ell'$.
These two situations have already been considered  in Case III and Case II
above (in that order): there we represented $z_\ell$ as $y_k \wedge y_m$
and $z_\ell'$ as $y_k' \wedge y_m'$. 
It follows 
from the inequality (\ref{changerelbdIII})
that the mass of the set of points $a$ that satisfy the first alternative
is at most $2 \gamma n < 112 \delta n$.
Similarly, from the inequality (\ref{changerelbdII})
and the fact that $\psi$ is measure-preserving,
the mass of the set of points $a$ that satisfy the second alternative
is also at most $2 \gamma n < 112 \delta n$.
Thus the total mass of the set of points of interest is at most $224 \delta n$. 
\end{proof}

Before completing the proof of Proposition~\ref{P:Feller}, we recall the
definition of the {\em Wasserstein metric}.  Suppose that $(E,d)$
is a complete, separable metric space.  Write $B$ for the set of continuous functions
functions $f:E \rightarrow \R$ such that $|f(x)| \le 1$ and $|f(x) - f(y)| \le d(x,y)$
for $x,y \in E$.  The Wasserstein (sometimes
transliterated as Vasershtein) distance between two Borel probability measures
$\alpha$ and $\beta$ on $E$ is given by
\be{defWass}
d_{\mathrm{W}} (\alpha,\beta) := \sup_{f \in B} \left|\int f d\alpha - \int f d\beta\right|.
\ee
The Wasserstein distance is a genuine metric on the space of Borel probability
measures and convergence with
respect to this distance implies weak convergence
(see, for example, Theorem 3.3.1 and Problem 3.11.2 of
\cite{EthierKurtz1986}).  If $V$ and $W$ are two 
$E$-valued random variables
on the same probability space $(\Sigma, {\mathcal{A}}, \mathbb{P})$ with 
distributions $\alpha$ and $\beta$, respectively, then
\be{Wasscompare}
\begin{split}
d_{\mathrm{W}} (\alpha,\beta) 
& \le \sup_{f \in B} \left|\mathbb{P}[f(V)] - \mathbb{P}[f(W)]\right| \\
& \le \sup_{f \in B} \mathbb{P}[|f(V)-f(W)|]
\le \mathbb{P}[d(V,W)]. \\
\end{split}
\ee

\noi
{\em Proof of Proposition~\ref{P:Feller}} \quad
For $(T,\rho) \in {\mathbf T}^{\mathrm{root}}$ and $t \ge 0$, let 
\be{defkernel}
{\mathbf{P}}_t((T,\rho), \boldsymbol{\cdot}) := {\bf P}^{(T,\rho)}\{X_t\in\boldsymbol{\cdot}\}.
\ee
We need to show that $(T,\rho) \mapsto {\mathbf{P}}_t((T,\rho), \boldsymbol{\cdot})$ is
weakly continuous for each $t \ge 0$.  This is equivalent to showing for each 
$(T,\rho) \in {\mathbf T}^{\mathrm{root}}$ and $t \ge 0$ that
\be{convinWass}
\lim_{(T',\rho') \rightarrow (T,\rho)} 
d_{\mathrm{W}} 
\left({\mathbf{P}}_t((T,\rho), \boldsymbol{\cdot}), {\mathbf{P}}_t((T',\rho'), \boldsymbol{\cdot})\right) = 0.
\ee

>From the coupling argument in the proof of part (iii) of Theorem~\ref{main}
(in particular, the inequality (\ref{uniformcontrol})), we have that
\be{Wasstriangle}
\begin{split}
& d_{\mathrm{W}} ({\mathbf{P}}_t((T,\rho), \boldsymbol{\cdot}), {\mathbf{P}}_t((T',\rho'), \boldsymbol{\cdot})) \\
& \quad \le
d_{\mathrm{W}} ({\mathbf{P}}_t((T,\rho), \boldsymbol{\cdot}), {\mathbf{P}}_t((R_\eta(T),\rho), \boldsymbol{\cdot})) \\
& \qquad + d_{\mathrm{W}} ({\mathbf{P}}_t((R_\eta(T),\rho), \boldsymbol{\cdot}), {\mathbf{P}}_t((R_\eta(T'),\rho'), \boldsymbol{\cdot})) \\
& \qquad + d_{\mathrm{W}} ({\mathbf{P}}_t((R_\eta(T'),\rho), \boldsymbol{\cdot}), {\mathbf{P}}_t((T',\rho'), \boldsymbol{\cdot})) \\
& \quad \le
d_{\mathrm{W}} ({\mathbf{P}}_t((R_\eta(T),\rho), \boldsymbol{\cdot}), {\mathbf{P}}_t((R_\eta(T'),\rho'), \boldsymbol{\cdot}))
+ 2 \eta.\\
\end{split}
\ee

By part (ii) of  Lemma~\ref{L1}, $R_\eta(T')$ converges to $R_\eta(T)$ 
as $(T',\rho')$ converges to $(T,\rho)$, and so it suffices to establish
(\ref{convinWass}) when $(T,\rho)$ and $(T',\rho')$
are finite trees, and so we will suppose this for the rest of the proof.

Fix $(T,\rho)$ and $\varepsilon>0$.  Suppose that $\delta > 0$
depending on $(T,\rho)$ and $\varepsilon$ is sufficiently small
that the conclusions of  Lemma~\ref{L4} hold
for any $(T',\rho')$ within distance $\delta$ of $(T,\rho)$.  
Let $(S,\rho)$ and $(S', \rho')$ be the subtrees guaranteed
by  Lemma~\ref{L4}.  From the coupling
argument in proof of part (iii) of Theorem~\ref{main}
we have 
\be{}
d_{\mathrm{W}} \left({\mathbf{P}}_t((T,\rho), \boldsymbol{\cdot}), {\mathbf{P}}_t((S,\rho), \boldsymbol{\cdot})\right) < \varepsilon
\ee
and
\be{}
d_{\mathrm{W}} \left({\mathbf{P}}_t((T',\rho'), \boldsymbol{\cdot}), {\mathbf{P}}_t((S',\rho'), \boldsymbol{\cdot}) \right)< \varepsilon.
\ee
It therefore suffices to give a bound on
$d_{\mathrm{W}} ({\mathbf{P}}_t((S,\rho), \boldsymbol{\cdot}), {\mathbf{P}}_t((S',\rho), \boldsymbol{\cdot}))$
that only depends on $\varepsilon$ and converges to zero as $\varepsilon$ converges to $0$.

Construct on some probability space $(\Sigma, {\mathcal{A}}, \mathbb{P})$ 
a Poisson point process $\Pi_0$ on the set $\R^{++} \times S^o$ 
with intensity $\lambda \otimes \mu$, where $\mu$ is the length
measure on $S$.  Construct on the same space another independent
Poisson point process on the set
$\{(t,x)\in\R^{++}\times\R^{++} : x \le t \}$ with
intensity $\lambda \otimes \lambda$ restricted to this set.
If we set 
$\Pi_0' := \{(t,\psi(x)) : (t,x) \in \Pi_0\} \subset \R^{++} \times (S')^o$,
then $\Pi_0'$ is a Poisson process on the set
$\R^{++} \times (S')^o$ 
with intensity $\lambda \otimes \mu'$, where $\mu'$ is the length
measure on $S'$ (because $\psi$ preserves length measure).
Now apply the construction of Subsection~\ref{SS:beginconstruct} to
realizations of $\Pi_0$ and $\Pi$ (respectively, $\Pi_0'$ and
$\Pi$) to get two ${\mathbf T}^{\mathrm{root}}$-valued processes
that we will denote by $(Y_t)_{t \ge 0}$ and $(Y_t')_{t \ge 0}$.
We see from the proof of Theorem~\ref{main} that $Y$ (respectively, $Y'$)
has the same law as $X$ under ${\bf P}^{(S,\rho)}$ 
(respectively, ${\bf P}^{(S',\rho')}$).

Define a map $\psi_t$ from $Y_t = S \amalg ]0,t]$ to $Y_t' = S' \amalg ]0,t]$
by setting the restriction of $\psi_t$ to $S$ be $\psi$ and the restriction
of $\psi_t$ to $]0,t]$ be the identity map.  
Let $d_t$ and $d_t'$ be the metrics on $Y_t$ and $Y_t'$, respectively.
We will bound the
rooted Gromov-Hausdorff distance between 
$Y_t$ and $Y_t'$ by bounding the distortion of $\psi_t$.

The cut-times for $Y$ and $Y'$ coincide.
If $\xi$ is a cut-point of $Y$ at some cut-time $\tau$,
then the corresponding cut-point
for $Y'$ will be $\psi(\xi)$.

It is clear that the distortion of $\psi_t$ is constant
between cut-times. 
Write $B_t$ for the set of points $b \in Y_t$
such that the subtree of $Y_t'$ above $\psi_t(b)$ is not  the image 
under $\psi_t$ of the subtree of $Y_t$ above $b$.
The set $B_t$ is unchanged between cut-times.

Consider a cut-time $\tau$ such that the corresponding cut-point
$\xi$ is in $Y_{\tau-} \setminus B_{\tau-}$.  If $x$ and $y$ are in the subtree
above $\xi$ in $Y_{\tau-}$, then they are moved together
by the re-grafting operation and their distance apart is unchanged
in $Y_\tau$.  Also, $\psi_{\tau-}(x)$ and $\psi_{\tau-}(y)$ are in 
subtree above $\psi_{\tau-}(\xi)$ in $Y_{\tau-}'$ 
and these two points are also moved together.
More precisely,
\be{}
d_\tau(x,y) = d_{\tau-}(x,y)
\ee
and
\be{}
d_\tau'(\psi_{\tau}(x),\psi_{\tau}(y)) = d_{\tau-}'(\psi_{\tau-}(x),\psi_{\tau-}(y)).
\ee
The same conclusion holds if neither 
$x$ or $y$ are in the subtree above $\xi$ in $Y_{\tau-}$.
If $x$ is in the subtree above $\xi$ in $Y_{\tau-}$ and
$y$ is not, then 
\be{}
d_\tau(x,y) = d_{\tau-}(x,\xi) + d_{\tau-}(\tau,y)
\ee
and
\be{}
d_\tau'(\psi_\tau(x), \psi_\tau(y)) = d_{\tau-}'(\psi_{\tau-}(x), \psi_{\tau-}(\xi)) + d_{\tau-}'(\tau, \psi_{\tau-}(y))
\ee
(where we recall that $\tau$ is the root in each of the trees
$Y_{\tau-}, Y_\tau, Y_{\tau-}', Y_\tau'$).
Combining these cases, we see that 
\be{disdoublesgood}
{\mathrm{dis}}(\psi_\tau) \le 2 {\mathrm{dis}}(\psi_{\tau-}).
\ee{}
Moreover, if $\xi \in Y_{\tau-} \setminus B_{\tau-}$, then
$B_{\tau} = B_{\tau-}$.

Also, for any $t \ge 0$ we always have the upper bound
\be{upperbdbad}
\begin{split}
{\mathrm{dis}}(\psi_t) 
& \le {\mathrm{diam}}(Y_t) + {\mathrm{diam}}(Y_t') \\
& \le {\mathrm{diam}}(S) + {\mathrm{diam}}(S') + 2t \\
& \le {\mathrm{diam}}(T) + {\mathrm{diam}}(T') + 2t \\
& \le 2 {\mathrm{diam}}(T) + d_{\mathrm{GH}^{\mathrm{root}}}((T,\rho), (T',\rho')) + 2t \\
& \le 2 {\mathrm{diam}}(T) + \delta + 2t \\
& =: D_t \\
\end{split}
\ee

Set $N_t := |\Pi_0 \cap (]0,t] \times S^o)| + |\Pi \cap \{(s,x) : 0 < x \le s \le t\}|$ 
and write $I_t$ for the indicator
of the event $\{\Pi_0 \cap (]0,t] \times B_0) \ne \emptyset\}$, which, by the above argument,
is the event that $\xi \in B_{\tau-}$ for some (cut-time, cut-point) pairs $(\tau,\xi)$ with $0 < \tau \le t$. 
We have
\be{}
\begin{split}
& d_{\mathrm{W}} ({\mathbf{P}}_t((S,\rho), \boldsymbol{\cdot}), {\mathbf{P}}_t((S',\rho), \boldsymbol{\cdot}) \\
& \quad \le \mathbb{P}[d_{\mathrm{GH}^{\mathrm{root}}}(Y_t,Y_t')] \\
& \quad \le \frac{1}{2} \mathbb{P}[{\mathrm{dis}}(\psi_t)] \\
& \quad \le \frac{1}{2} \mathbb{P}\left[\varepsilon 2^{N_t}  + I_t D_t\right]\\
& \quad = \frac{1}{2}\left\{\varepsilon \exp\left(\mu(T) t + \frac{t^2}{2}\right) 
+ \left[1 - \exp\left( - \varepsilon t \right)\right] D_t\right\},\\
\end{split}
\ee
and this suffices to complete the proof.
\hfill $\qed$

\section{Asymptotics of the Aldous-Broder algorithm}
\label{S:Aldous-Broder}

Given an irreducible Markov matrix $\mathbb{P}$ with 
state space $V$, there is a natural probability measure
on the collection of combinatorial trees with vertices
labeled by $V$ that assigns mass
\be{wform}
   C^{-1}\prod \mathbb{P}(x,y)
\ee
to the tree $T$, where $C$ is a normalization constant
and the product is over pairs of adjacent
vertices $(x,y)$ in $T$ ordered so that $y$ is on the path
from the root to $x$.
For example, if
$\mathbb{P}(x,y)\equiv 1/|V|$ 
for all $x,y\in V$, (so that the associated Markov chain
consists of successive uniform random picks from $V$), then
the distribution (\ref{wform}) is uniform on the set of
$|V|^{|V|-1}$ rooted combinatorial trees labeled by $V$.

The {\em Aldous-Broder algorithm} \cite{AT88,bro89,MR91h:60013}
is a tree-valued Markov chain that has the distribution
in (\ref{wform}) as its stationary distribution.
The discrete time version of the algorithm has the
following transition dynamics.
\begin{itemize}
\item
Pick a vertex $\upsilon$ at random according to
${\mathbb{P}}(\rho, \cdot)$, where $\rho$ is the current root.
\item
If $\upsilon = \rho$, do nothing.
\item
If $\upsilon \ne \rho$:
\begin{itemize}
\item
Erase the edge connecting $\upsilon$ to the unique vertex 
adjacent to
$\upsilon$ and on the path from $\rho$ to $\upsilon$.
\item
Insert a new edge between $\upsilon$ and $\rho$.
\item
Designate $\upsilon$ as the new root.
\end{itemize}
\end{itemize}

\setlength{\unitlength}{0.7pt}
\begin{picture}(300,380)(-20,-300)
\thinlines
\put(0,-15){\line(1,0){90}}
\put(-18,46){\mbox{$\rho_1$}}
\multiput(90,15)(-18,6){5}{\line(-3,1){15}}
\put(0,45){\line(1,0){90}}
\put(0,45){\line(0,-1){30}}
\put(90,15){\line(0,-1){30}}
\put(90,15){\circle{3}}
\put(120,15){\circle{3}}
\put(194,25){\mbox{$\rho_2$}}
\put(120,-15){\line(1,0){90}}
\put(120,45){\line(1,0){90}}
\multiput(120,45)(0,-4){8}{\line(0,-1){2}}
\put(210,15){\line(0,-1){30}}
\put(210,15){\line(-3,1){90}}
\put(330,-15){\circle{3}}
\put(250,24){\mbox{$\rho_3$}}
\put(240,-15){\line(1,0){90}}
\put(240,45){\line(1,0){90}}
\multiput(330,15)(0,-4){8}{\line(0,-1){2}}
\put(330,15){\line(-3,1){90}}
\put(440,0){\mbox{$\rho_4$}}
\put(360,-15){\line(1,0){90}}
\put(360,45){\line(1,0){90}}
\put(360,15){\line(1,0){90}}
\put(450,15){\line(-3,1){90}}
\put(30,-120){\circle{3}}
\put(49,-153){\mbox{$\rho_1$}}
\put(30,-120){\line(0,1){60}}
\multiput(30,-150)(0,4){8}{\line(0,1){2}}
\put(30,-150){\line(1,1){20}}
\put(30,-150){\line(-1,1){20}}
\put(160,-90){\circle{3}}
\put(175,-128){\mbox{$\rho_1$}}
\put(175,-153){\mbox{$\rho_2$}}
\put(160,-150){\line(-1,1){40}}
\put(160,-120){\line(1,1){20}}
\multiput(160,-120)(0,4){8}{\line(0,1){2}}
\put(290,-90){\circle{3}}
\put(320,-98){\mbox{$\rho_1$}}
\put(300,-123){\mbox{$\rho_2$}}
\put(300,-153){\mbox{$\rho_3$}}
\put(290,-120){\line(1,1){40}}
\put(290,-90){\line(0,1){30}}
\multiput(290,-120)(0,4){8}{\line(0,1){2}}
%
\put(435,-68){\mbox{$\rho_1$}}
\put(435,-93){\mbox{$\rho_2$}}
\put(435,-123){\mbox{$\rho_3$}}
\put(435,-153){\mbox{$\rho_4$}}
\put(420,-120){\line(0,1){60}}
\put(420,-60){\line(1,1){20}}
\put(420,-150){\line(-1,1){20}}
\thicklines 
\put(210,15){\line(-3,1){90}}
\put(290,-150){\line(0,1){30}}
\put(160,-150){\line(0,1){30}}
\put(420,-150){\line(0,1){30}}
\put(240,15){\line(1,0){90}}
\put(360,15){\line(3,-1){90}}
\put(-20,-190){\makebox{
\begin{minipage}[t]{12cm}
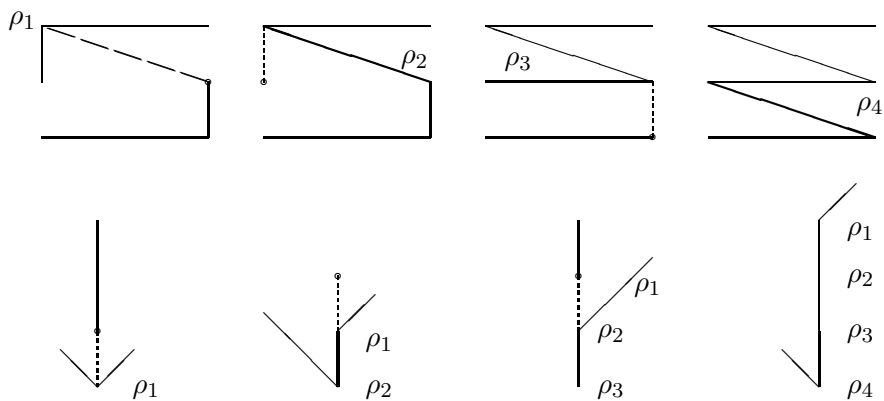
{\em {\sc Figure~4} illustrates  
Aldous-Broder algorithm 
in discrete time. Here the dots mark the vertices that become the root
in the next step. Once the new root is chosen
we erase the edge adjacent to the new root on the path 
from the old to the new root and insert an edge connecting the
old to the new root. The bold edges in the picture present the 
edges introduced recently.}\end{minipage}}}
\end{picture}

It will be more convenient for us to work with the continuous
time version of this algorithm in which the above transitions
are made at the arrival times of an independent 
Poisson process with rate $|V|/(|V|-1)$ (so that
the continuous time chain makes actual jumps at rate $1$).

We can associate a rooted compact real tree with
a rooted labeled combinatorial tree in the
obvious way by
thinking of the edges as line segments with length $1$.
Because we don't record the labeling, the process that arises
from mapping the continuous-time Aldous-Broder algorithm
in this way won't be Markovian in general.  However, this process
will be Markovian in the case where $\mathbb{P}$
is the transition matrix for i.i.d. uniform sampling
(that is, when $\mathbb{P}(x,y) = 1/|V|$ for all
$x,y \in V$) and we assume this from now on.  
The following result says that if
we rescale ``space'' and time appropriately, then
this process converges to the root growth
with re-grafting process.  If $T = (T,d,\rho)$
is a rooted compact real tree and $c>0$, we
write $c T$ for the tree $(T,c\,d,\rho)$ (that is,
$c T=T$ as sets and the roots are the same, but
the metric is re-scaled by $c$).

\begin{proposition}
\label{ABtoRGR}
Let $Y^n = (Y_t^n)_{t \ge 0}$ be a sequence of Markov processes
that take values in the space of rooted compact real trees with
integer edge lengths and evolve according to the dynamics
associated with the continuous-time Aldous-Broder chain
for i.i.d. uniform sampling.  Suppose that each tree
$Y_0^n$ is non-random with total branch length 
$N_n$, that $N_n$ converges to infinity
as $n \rightarrow \infty$,
and that $N_n^{-1/2} Y_0^n$ converges in the
rooted Gromov-Hausdorff metric to some rooted compact
real tree $T$ as $n \rightarrow \infty$.  
Then, in the sense of weak convergence
of processes on the space of c\`adl\`ag paths equipped
with the Skorohod topology, $(N_n^{-1/2} Y^n(N_n^{1/2}t))_{t \ge 0}$
converges as $n \rightarrow \infty$ to the root growth with re-grafting process
$X$ under ${\mathbf{P}}^T$.
\end{proposition}

\begin{proof}
Define $Z^n = (Z_t^n)_{t \ge 0}$ by
\be{}
Z_t^n := N_n^{-1/2} Y^n(N_n^{1/2}t).
\ee
For $\eta>0$, let $Z^{\eta,n}$ be the ${\mathbf{T}}^{\mathrm{root}}$-valued
process constructed as follows.
\begin{itemize}
\item
Set $Z_0^{\eta,n} = R_{\eta_n} (Z_0^n)$,
where $\eta_n := N_n^{-1/2}\lfloor N_n^{1/2} \eta \rfloor$.
\item
The value of $Z^{\eta,n}$ is unchanged between jump times of
$(Z_t^n)_{t \ge 0}$.
\item
At a jump time $\tau$ for $(Z_t^n)_{t \ge 0}$,
the tree  $Z_\tau^{\eta,n}$ is the subtree of
$Z_\tau^n$ spanned by $Z_{\tau-}^{\eta,n}$
and the root of $Z_\tau^n$.
\end{itemize}
An argument similar to that in the proof of Lemma \ref{L0}
shows that
\be{}
\sup_{t \ge 0} d_{\mathrm{H}}(Z_t^n, Z_t^{\eta,n}) \le \eta_n,
\ee
and so it suffices to show that
$Z^{\eta,n}$ converges weakly as $n \rightarrow \infty$ to $X$ under
${\mathbf{P}}^{R_\eta(T)}$.

Note that 
$Z_0^{\eta,n}$ converges to $R_\eta(T)$
as $n \rightarrow \infty$.  Moreover, 
if $\Lambda$ is the map that sends a tree to its total length 
(that is, the total
mass of its length measure),
then 
$\lim_{n \rightarrow \infty} \Lambda(Z_0^{\eta,n}) 
= \Lambda \circ R_\eta(T) < \infty$
by Lemma \ref{contlength} below.

The pure jump process $Z^{\eta,n}$ is clearly Markovian.
If it is in a state $(T',\rho')$, then it jumps with the following
rates.
\begin{itemize}
\item
With rate $N_n^{1/2} (N_n^{1/2} \Lambda(T'))/N_n = \Lambda(T')$,
one of the $N_n^{1/2} \Lambda(T')$ points in $T'$ that are at
distance a positive integer multiple of $N_n^{-1/2}$ from the root $\rho'$
is chosen uniformly at random and the subtree above
this point is joined to $\rho'$ by an edge of length
$N_n^{-1/2}$.  The chosen point becomes the new root
and an arc of length  $N_n^{-1/2}$
that previously led from the new root toward
$\rho'$ is erased.  Such a transition
results in a tree with the same total length as $T'$.
\item
With rate  $N_n^{1/2} - \Lambda(T')$, a new root not present
in $T'$ is attached to $\rho'$ by an edge of length 
$N_n^{-1/2}$.  This results in a tree with total length
$\Lambda(T') + N_n^{-1/2}$.
\end{itemize}
It is clear that these dynamics converge to those
of the root growth with re-grafting process, 
with the first class of transitions
leading to re-graftings in the limit and the second class
leading to root growth.
\end{proof}

\begin{lemma}
\label{lengthformula}
Let $(T, d, \rho) \in {\mathbf T}^{\mathrm{root}}$ and suppose that 
$\{x_0, \ldots, x_n\} \subset T$ spans $T$, so that the root $\rho$
and the leaves of $T$ form a subset of $\{x_0, \ldots, x_n\}$.  Then the total
length of $T$ (that is, the total mass of its length measure)
is given by
\begin{equation*}
d(x_0, x_1) +
\sum_{k=2}^n \bigwedge_{0 \le i < j \le k-1} \frac{1}{2}\left(d(x_k,x_i) + d(x_k, x_j) - d(x_i, x_j)\right).
\end{equation*}
\end{lemma}

\begin{proof}
This follows from the observation that 
the distance from the point $x_k$
to the arc $[x_i, x_j]$ is
\begin{equation*}
\frac{1}{2}\left(d(x_k,x_i) + d(x_k, x_j) - d(x_i, x_j)\right),
\end{equation*}
and so length of the arc connecting $x_k$, $2 \le k \le n$, to the subtree
spanned by $x_0, \ldots, x_{k-1}$ is
\begin{equation*}
\bigwedge_{0 \le i < j \le k-1} \frac{1}{2}\left(d(x_k,x_i) + d(x_k, x_j) - d(x_i, x_j)\right).
\end{equation*}
\end{proof}

\begin{lemma}
\label{contlength}
Let $\Lambda: {\mathbf T}^{\mathrm{root}} \rightarrow \R \cup \{\infty\}$
be the map that sends a tree to its total length.  For $\eta>0$, the map $\Lambda \circ R_\eta$
is continuous.
\end{lemma}

\begin{proof}
For all $\eta>0$ we have by Lemma~\ref{L1} that:
$R_\eta = R_{\eta/2} \circ R_{\eta/2}$, the map $R_\eta$ is continuous,
and the range of $R_\eta$ consists of finite trees.
It therefore suffices to show for all $\eta>0$ that if 
$(T, d, \rho)$ is a fixed finite tree
and $(T', d', \rho')$ is any another finite tree sufficiently close to $T$,
then $\Lambda \circ R_\eta(T')$ is close to $\Lambda \circ R_\eta(T)$.

Suppose, therefore, that $(T, d, \rho)$ is a fixed finite tree 
with leaves $\{x_1, \ldots, x_n\}$ and that
$(T', d', \rho')$ is another finite tree with 
\begin{equation*}
d_{{\mathrm{GH}}^{\mathrm{root}}}((T, d, \rho),(T', d', \rho'))<\delta,
\end{equation*}
where $\delta$ is small enough that the conclusions
of Lemma~\ref{L:Tdoubprimeexist} hold.  Consider a rooted subtree 
$(T'', d', \rho')$ 
of $(T', d', \rho')$ and a map $\bar f: T \rightarrow T''$ with the
properties guaranteed by Lemma~\ref{L:Tdoubprimeexist}.
Set $x_k' = \bar f(x_k)$ for $1 \le k \le n$.

Fix $\kappa>0$.
For $1 \le k \le n$, write $\hat x_k \in T$ for the point on the arc 
$[\rho, x_k]$ that is at distance $\kappa \wedge d(\rho, x_k)$ from
$x_k$.  Set $\hat x_0 := \rho$.  Define $\hat x_0', \ldots, \hat x_n' \in T''$
similarly.  Note that $R_\kappa(T)$ is spanned by $\{\hat x_0, \ldots, \hat x_n\}$
and $R_\kappa(T'')$ is spanned by $\{\hat x_0', \ldots, \hat x_n'\}$.
By Lemma~\ref{lengthformula},
\begin{equation*}
\Lambda \circ R_\kappa(T)
=
d(\hat x_0, \hat x_1) +
\sum_{k=2}^n \bigwedge_{0 \le i < j \le k-1} \frac{1}{2}\left(d(\hat x_k,\hat x_i) + d(\hat x_k, \hat x_j) - d(\hat x_i, \hat x_j)\right).
\end{equation*}
and
\begin{equation*}
\Lambda \circ R_\kappa(T'')
=
d'(\hat x_0', \hat x_1') +
\sum_{k=2}^n \bigwedge_{0 \le i < j \le k-1} \frac{1}{2}\left(d'(\hat x_k',\hat x_i') + d'(\hat x_k', \hat x_j') - d'(\hat x_i', \hat x_j')\right).
\end{equation*}
Also observe that
\begin{equation*}
\begin{split}
d(\hat x_i, \hat x_j) & = (d(x_0,  x_i) - \kappa)_+ + (d(x_0, x_j) - \kappa)_+ \\
& \quad - 2 \biggl[(d(x_0,  x_i) - \kappa)_+ \wedge (d( x_0,  x_j) - \kappa)_+  \\
& \qquad \wedge \Bigl\{\frac{1}{2}\bigl(d( x_0,  x_i) + d(x_0,  x_j) - d(x_i,  x_j)\bigr)\Bigr\}\biggr] \\
\end{split}
\end{equation*}
and
\begin{equation*}
\begin{split}
d'(\hat x_i',  \hat x_j') & = (d'( x_0',  x_i') - \kappa)_+ + (d'( x_0',  x_j') - \kappa)_+ \\
& \quad - 2 \biggl[(d'( x_0',  x_i') - \kappa)_+ \wedge (d'( x_0',  x_j') - \kappa)_+  \\
&  \qquad \wedge \Bigl\{\frac{1}{2}\bigl(d'( x_0',  x_i') + d'( x_0',  x_j') - d'( x_i',  x_j')\bigr)\Bigr\}\biggr].\\
\end{split}
\end{equation*}

Now the function $t \mapsto (t - \kappa)_+$ is Lipschitz with Lipschitz constant $1$ for all $\kappa>0$,
and it follows that there is a family of Lipschitz functions $F_\kappa$, $\kappa>0$, with 
Lipschitz constants uniformly bounded by some constant $C$ such that
\begin{equation*}
\Lambda \circ R_\kappa(T)
= F_\kappa\left(\left(d(x_i, x_j)\right)_{0 \le i,j \le n}\right)
\end{equation*}
and
\begin{equation*}
\Lambda \circ R_\kappa(T'')
= F_\kappa\left(\left(d'(x_i', x_j')\right)_{0 \le i,j \le n}\right).
\end{equation*}
By construction $|d(x_i, x_j) - d'(x_i', x_j')| < 8 \delta$, and so
\begin{equation*}
|\Lambda \circ R_\kappa(T) - \Lambda \circ R_\kappa(T'')| \le 8\delta C
\end{equation*}
for all $\kappa > 0$.

Because $d_{\mathrm{H}}(T',T'') < 3 \delta$, we have
\begin{equation*}
\Lambda \circ R_\eta(T'') 
\le \Lambda \circ R_\eta(T') 
\le \Lambda \circ R_{\eta - 3 \delta}(T'').
\end{equation*}
Thus
\begin{equation*}
\Lambda \circ R_\eta(T) - 8 \delta C
\le \Lambda \circ R_\eta(T') 
\le \Lambda \circ R_{\eta - 3 \delta}(T) + 8 \delta C.
\end{equation*}
Since
$\lim_{\delta \downarrow 0} \Lambda \circ R_{\eta - 3 \delta}(T)
= \Lambda \circ R_\eta(T)$, this suffices to establish the result.
\end{proof}

An alternative algorithm for simulating from the distribution in (\ref{wform})
in the case of i.i.d. uniform sampling
is the complete graph special case
of Wilson's loop-erased walk algorithm for generating a uniform
spanning tree of a graph \cite{MR99g:60116, MR1427525, MR96m:68144}.  Asymptotics of
the latter algorithm have been investigated in \cite{MR2003j:60011}.
Wilson's algorithm was also used in \cite{math.PR/0410430} to show
that the finite-dimensional distributions of the re-scaled
uniform random spanning tree for the $d$-dimensional discrete
torus converges to the Brownian CRT as the number of vertices
goes to infinity when $d \ge 5$.

\section{Rayleigh process}
\label{S:Rayleigh}

Suppose that we take the root growth with re-grafting process
$(X_t)_{t \ge 0}$ under ${\mathbf{P}^T}$ for some 
$T \in {\mathbf{T}}^{\mathrm{root}}$, we fix a point
$x \in T$, and we denote by $R_t$ the distance between $x$
and the root $t$ of $X_t$ (that is, $R_t$ is the height of $x$
in $X_t$).  According to the root growth with
re-grafting dynamics,
$R_t$ grows deterministically with unit
speed between cut-time $\tau$ for which the
corresponding cut-point falls on the arc $[\tau,x]$.
Such cut-times $\tau$ come along at intensity $R_{t-} \, dt$
in time, and at $\tau-$ the position of the corresponding cut-point
is uniformly distributed on the arc $[\tau, x]$
conditional on the past up to $\tau-$, so that $R_{\tau}$ is uniformly
distributed on $[0, R_{\tau-}]$ conditional on the past up to $\tau-$.
Consequently, the $\R^+$-valued process $(R_t)_{t \ge 0}$
is autonomously Markovian.  In particular, 
$(R_t)_{t \ge 0}$ is an example of the class of
{\em piecewise deterministic} Markov processes discussed
in the Introduction.

In order to describe the properties of $(R_t)_{t \ge 0}$, we need the following
definitions.
A non-negative random variable $R$ is
said to have {\em standard Rayleigh distribution} if
it is distributed as the length of a standard normal vector in $\R^2$, 
that is,
\be{Rdef}
   {\mathbb P}\{R>r\} 
 = 
   \exp\left(-\frac{r^2}{2}\right), \quad r\ge 0.
\ee
If $R^*$ is distributed according to the size-biased standard Rayleigh distribution,
that is, 
\be{Rast}
   {\mathbb P}\{R^\ast\in\mathrm{d}r\}
=\frac{r{\mathbb P}\{R\in\mathrm{d}r\}}{{\mathbb
     P}[R]}=\sqrt{\frac{2}{\pi}}\,r^2e^{-\frac12 r^2}\mathrm{d}r, \quad r\ge 0,
\ee
and if $U$ is a uniform random variable that is independent of
$R^\ast$,
then $UR^\ast$ has the inverse size-biased standard Rayleigh distribution:
\be{URast}
   {\mathbb P}\{U R^\ast\in\mathrm{d}r\}=\frac{r^{-1}{\mathbb P}\{R\in\mathrm{d}r\}}
   {{\mathbb P}[R^{-1}]}=\sqrt{\frac{2}{\pi}}\,e^{-\frac12 r^2}
    \mathrm{d}r, \quad r\ge 0.
\ee
Thus $R^{\ast}$ and $U R^\ast$ are distributed as the length of a
standard normal vector in $\R^3$ and $\R$, respectively.
\bi

For reasons that are apparent from Proposition~\ref{rayproc}
below, we call the process $(R_t)_{t\ge 0}$ the
{\em Rayleigh process}. We note that there is
a body of literature on stationary processes
with Rayleigh one-dimensional marginal distributions
that arise as the length process of a vector-valued process in $\R^2$ with
coordinate processes that are independent copies of
some stationary centered Gaussian process (see, for example,
\cite{MR45:7804,MR20:1371,MR2003g:60001}).

\begin{proposition}\label{rayproc}
Consider the Rayleigh process $(R_t)_{t\ge 0}$.
Write ${\mathbf P}^r$ for the law of $(R_t)_{t\ge 0}$
started at $r \ge 0$.
\begin{itemize}
\item[(i)]
The unique stationary distribution of the Rayleigh process
is the standard Rayleigh distribution and
the total variation distance between ${\bf P}^r\{R_t \in \cdot\}$ and 
the standard Rayleigh distribution
converges to $0$ as $t \rightarrow \infty$.
\item[(ii)] Under ${\mathbf P}^0$,
for each fixed $t>0$, $R_t$ has the same law as $R\wedge t$,
where $R$ has the standard Rayleigh distribution.
\item[(iii)]
For $x>0$, the mean return time to $x$ is 
$x^{-1}e^{\frac{1}{2}x^2}$.
\item[(iv)] If $\tau_n$ denotes the $n$th jump time of $(R_t)_{t\ge
  0}$, then as $n \rightarrow \infty$ the triple
$(R_{\tau_n},R_{\tau_{n+1}-},R_{\tau_{n+1}})$
converges in law to the triple
$ (U' R^\ast,R^\ast,U'' R^\ast)$,
where $U'$ and $U''$ are 
independent uniform random variables on $]0,1[$
independent of $R^\ast$, and $R^\ast$ has
the size-biased Rayleigh distribution.
\item[(v)] The jump counting process $N(t):= |\{n \in \N : \tau_n \le t\}|$
has asymptotically stationary increments under
${\mathbf P}^r$ for any $r \ge 0$, and 
\be{count}
   \frac{1}{t}N(t)\rightarrow\sqrt{\frac{\pi}{2}},
\quad {\mathbf P}^r-{\mathrm{a.s.}}
\ee
as $t \rightarrow \infty$.
\end{itemize}
\end{proposition}

\begin{proof} 
(i) Let $\bar \Pi$ be a Poisson point process in $\R \times \R_+$
with Lebesgue intensity.  For $-\infty<t<\infty$ let
\begin{equation}
\bar R_t := \inf\{x + (t-s) : (s,x) \in \bar \Pi, \, s \le t\}.
\end{equation}
It is clear that $(\bar R_t)_{t \in \R}$ is a stationary Markov process
with the transition dynamics of the Rayleigh process.
Similarly, for $r \in \R_+$ and $t \ge 0$, set 
\begin{equation}
R_t^r = (r + t) \wedge \inf\{x + (t-s) : (s,x) \in \bar \Pi, \, 0 \le s \le t\}.
\end{equation}
Then $(R_t^r)_{t \ge 0}$ has the same law as the Rayleigh process
under ${\bf P}^r$.

Note that the event $\{\bar R_t > r\}$ is the event that $\bar \Pi$
has no points in the triangle with vertices $(t-r,0), (t,0), (t,r)$
and area $r^2/2$. Thus ${\mathbf P}\{\bar R_t > r\} = \exp(-r^2/2)$
and the standard Rayleigh distribution is a stationary distribution for the
Rayleigh process.

Let $T^r := \inf\{t \ge 0 : R_t^r = \bar R_t\}$.  Note that 
$R_t^r = \bar R_t$ for all $t \ge T^r$.  
Note that $T^r > t$ if and only if either $\bar R_0 > r$
and $\bar \Pi$ puts no points into the
quadrilateral with vertices $(0,0)$, $(t, 0)$, $(t, r+t)$, $(0,r)$,
or $\bar R_0 \le r$ and
$\bar \Pi$ puts no points into the
quadrilateral with vertices $(0,0)$, $(t, 0)$, $(t, \bar R_0+t)$, 
$(0,\bar R_0)$.
Hence
\begin{equation}
\begin{split}
{\mathbb P}\{T^r > t\} & = \exp\left(-\frac{r^2}{2}\right) \exp\left(- \frac{1}{2}(r + (r+t))t\right) \\
& \quad + \int_0^r \exp\left(- \frac{1}{2}(x + (x+t))t\right) \, x \exp\left(-\frac{x^2}{2}\right) \, d x. \\
\end{split}
\end{equation}
By the standard coupling inequality, the total variation between
${\mathbf P}^r\{R_t \in \cdot\}$ and ${\mathbb P}\{\bar R_t \in \cdot\}$
is at most $2 {\mathbb P}\{T^r>t\}$, which converges to $0$ as $t \rightarrow \infty$.
This certainly shows that the standard Rayleigh
distribution is the unique stationary distribution.

\noi
(ii) Note that $R_t^r>x$ if and only if $r+t \ge t > x$ and
there are no points of $\bar \Pi$ in the triangle
with vertices $(t-x,0), (t,0), (t,x)$ of area $x^2/2$
or $r+t>x \ge t$ and there are not points of $\bar \Pi$
in the  quadrilateral with vertices $(0,0), (t,0), (t,x), (0,x-t)$
of area $((x-t) + x)t/2 = x^2/2 - (x-t)^2/2$.
In either case,
\be{hazr}
   {\mathbf P}^{r}\{R_t>x\}
 =
   1\{r+t>x\}\exp\left(-\frac{1}{2}x^2+\frac{1}{2}((x-t)_+)^2\right).
\ee
Taking $r=0$ gives the result.\sm

\noi
(iii) Let $T_y:=\inf\{t>0:\,R_t=y\}$.
It is obvious from the Poisson construction that, for all $x \ge 0$ and $y>0$,
${\mathbf P}^x\{0<T_y<\infty\}=1$ and
${\mathbf P}^x[\exp(u T_y)] < \infty$
for all $u$ in some neighbourhood of $0$.
The Laplace transforms ${\mathbf P}^x[\exp-\lambda T_y]$ are determined by
standard methods of renewal theory: 
\be{lts}
   {\mathbf P}^x[\exp-\lambda T_x]=\frac{U_x(\lambda)}{1 + U_x(\lambda)},
   \quad \lambda >0,
\ee
where by (\ref{hazr}),
\be{Ux}
\begin{aligned}
   U_x(\lambda)
 &:=
   \int_0^\infty\mathrm{d}t\,e^{-\lambda t } \frac{{\mathbf P}^x\{R_t \in dx\}}{dx}
  \\
 &=
   \int_0^x\mathrm{d}t\,e^{-\lambda t }\,t\,e^{-xt+t^2/2}+
   {\mathbb P}\{R\in \mathrm{d}x\}\frac{e^{-\lambda x}}{\lambda }.
\end{aligned}
\ee
In particular, it follows easily that the mean return time of state
$x$
\be{means}
   {\mathbf P}^x[T_x]= 
   -\lim_{\lambda \downarrow 0}\frac{1}{\lambda}{\mathbf P}^x[\exp-\lambda T_x]
\ee
is the inverse of the density of $R$ at $x$, that is
$x^{-1} e^{\frac{1}{2}x^2}$, as claimed.
\sm

\noi
(iv) Let $\bar \tau := \inf\{t>0 : \bar R_t \ne \bar R_{t-}\}$.
By part (i), the joint distribution of $(R_{\tau_n},R_{\tau_{n+1}-},R_{\tau_{n+1}})$ converges
to the joint distribution of 
$(\bar R_0, \bar R_{\bar \tau -}, \bar R_{\bar \tau})$ conditional
on $\bar R_0 \ne \bar R_{0-}$. Let $C$ denote the intensity of the
stationary point process $\{t \in \R : \bar R_t \ne \bar R_{t-}\}$.
Then
\begin{equation}
\begin{split}
&{\mathbb P}\{\bar R_0 \in dx, \bar R_{\bar \tau -} \in dy, \bar R_{\bar \tau} \in dz \, | \, \bar R_0 \ne \bar R_{0-}\} \\
& \quad =
C^{-1} \exp\left(-\frac{1}{2} x^2\right) \, dx 
\, \exp\left(-\frac{1}{2}(x+y)(y-x)\right) \, dy \, dz \\
& \quad =
 \left[ \frac{1}{y} \, dx \right]
 \times \left[C^{-1} y^2 \exp\left(-\frac{1}{2} y^2\right) \, dy \right]
 \times \left[\frac{1}{y} \, dz \right]\\
\end{split}
\end{equation}
for $x < y$ and $z < y$.
The result now follows from (\ref{Rast}), which also identifies
$C = \sqrt{\pi/2}$.

\noi
(v) The stationary point process $\{t \in \R : \bar R_t \ne \bar R_{t-}\}$
is clearly ergodic by construction, and it has intensity $\sqrt{\pi/2}$
from the argument in part (iv).  For any $r>0$, it follows from the
argument in part (i) that $R_t^r = \bar R_t$ for all $t$ sufficiently
large, and so the result follows from the ergodic theorem applied to
$\{t \in \R : \bar R_t \ne \bar R_{t-}\}$.
\end{proof}

Then the following corollary is a
consequence of Proposition~\ref{ABtoRGR}.
See \cite{MR2003f:60168}, where similar scaling limits are derived.

\begin{cor}\label{prpscale}
For each $N\in\N$, let $(\tilde R_t^N)_{t \ge 0}$ denote
a continuous time
Markov chain with state space $\{1,\ldots,N\}$
and infinitesimal generator matrix
\be{Q}
   \tilde Q^N(i,j)
 :=
\begin{cases}
  1/N, & 1\le j\le i-1,\\ 
   -(N-1)/N, & j=i,\\
   (N-i)/N, & j=i+1,\\
   0, & \text{otherwise.}
\end{cases}.
\ee
Write $\tilde{\mathbf P}^{N,r}$, $r \in  \{1,\ldots,N\}$,
for the corresponding family of laws.
If a sequence $(r_N)_{N \in \N}$, $r_N \in  \{1,\ldots,N\}$,
is such that $\lim_{N \rightarrow \infty} N^{-1/2} r_N = r_\infty$
exists, then the law of 
$N^{-1/2}\left(\tilde R_{t \sqrt{N}}^N\right)_{t \ge 0}$
under $\tilde{\mathbf P}^{N,r_N}$ converges to that of
the Rayleigh process $(R_t)_{t\ge 0}$ under ${\mathbf P}^{r_\infty}$
in the usual sense of convergence  of
c\`adl\`ag  processes with the Skorohod topology.
\end{cor}

\bigskip
\noindent
{\bf Acknowledgment:}  We thank Jean-Fran\c{c}ois Le Gall and an anonymous referee for
several helpful comments that improved the presentation and correctness
of the paper.


\end{document}